\documentclass[11pt]{amsart}

\usepackage{amsaddr}
\usepackage{graphicx}
\usepackage[margin=1in]{geometry}
\usepackage[font=footnotesize, labelfont=bf]{caption}
\usepackage{pgfplots}
\usepackage{mathtools}
\usepackage[inline]{enumitem}
\usepackage{siunitx}
\usepackage{booktabs}
\usepackage{algorithm}
\usepackage{algpseudocode}

\usepackage{etoolbox}
\AtBeginEnvironment{algorithmic}{\small}

\usepackage[table]{xcolor}
\usepackage[sortcites=true,sorting=none,
citestyle=numeric-comp,backref=true,maxbibnames=10]{biblatex}
\definecolor{blue}{RGB}{49,114,174}
\definecolor{red}{RGB}{200,16,46}
\usepackage[colorlinks=true,
linkcolor=blue,
citecolor=blue]{hyperref}

\usepackage{cleveref}

\crefname{equation}{}{}
\Crefname{equation}{}{}

\newcommand\fnum[1]{\num[scientific-notation=fixed,round-precision=3,round-mode=places]{#1}}
\newcommand\snum[1]{\num[round-precision=3,round-mode=places,tight-spacing=true,retain-zero-exponent=true,retain-explicit-plus=true,output-exponent-marker=\text{e}]{#1}}

\newcommand\inum[1]{\num[scientific-notation=fixed,round-precision=0,round-mode=places,group-separator = {,}]{#1}}

\makeatletter
\renewcommand{\subsection}{%
  \@startsection{subsection}{2}%
  {\z@}%
  {1.5ex \@plus .2ex}%
  {.8ex}%
  {\normalfont\bfseries\itshape}}

\renewcommand{\subsubsection}{%
  \@startsection{subsubsection}{3}%
  {\z@}%
  {1.25ex \@plus .2ex}%
  {-1em}%
  {\normalfont\bfseries\itshape}}

\renewcommand{\paragraph}{%
  \@startsection{paragraph}{4}%
  {\z@}%
  {1ex}%
  {-1em}%
  {\normalfont\bfseries\itshape}}
\makeatother

\newcommand{\defeq}{\ensuremath{\mathrel{\mathop:}=}}

\newcommand\ui{\ensuremath{u_{\textit{i}}}} 
\newcommand\up{\ensuremath{u_{\textit{p}}}} 
\newcommand\un{\ensuremath{u_{\textit{n}}}} 
\newcommand\uo{\ensuremath{u_{\textit{o}}}} 
\newcommand\uh{\ensuremath{u_{\textit{h}}}} 
\newcommand\ug{\ensuremath{u_{\textit{g}}}} 
\newcommand\uw{\ensuremath{u_{\textit{w}}}} 
\newcommand\ipoint[1]{\noindent\textbf{#1}.}

\newcommand\igrad{\ensuremath{\nabla}}
\newcommand{\idiv}{\ensuremath{\nabla\cdot}}

\newcommand{\rom}{ROM}
\newcommand{\trom}{TROM}
\newcommand{\htrom}{HOSVD-TROM}
\newcommand{\ttrom}{TT-ROM}
\newcommand{\podrom}{POD-ROM}
\newcommand{\fom}{FOM}

\DeclareMathOperator*{\minopt}{minimize}
\newcommand\tabadjust{\centering\footnotesize}

\title[Tensorial ROM for Brain Tumor Growth Modeling]{Tensorial Reduced-Order Models for Parametric Coupled Reaction-Diffusion Systems: Application to Brain Tumor Growth Modeling}
\author{Asikul Islam, Md Rezwan Bin Mizan, Maxim Olshanskii \and Andreas Mang}
\address{Department of Mathematics, University of Houston}

\email{\{molshan, andreas\}@math.uh.edu}
\date{\today}

\begin{document}

\maketitle

\begin{abstract}
We construct efficient surrogate models for parametric forward operators arising in brain tumor growth simulations, governed by coupled semilinear parabolic reaction--diffusion systems on heterogeneous two- and three-dimensional domains. We consider two models of increasing complexity: a scalar single-species formulation and a six-state, nine-parameter multi-species go-or-grow model. The governing equations are discretized using a finite volume method and integrated in time via an operator-splitting strategy. We develop tensorial reduced-order model (\trom) surrogates based on the Higher-Order Singular Value Decomposition in Tucker format and the Tensor Train decomposition, each in intrusive and non-intrusive variants. The models are compared against a classical proper orthogonal decomposition (POD) ROM baseline. Numerical experiments with up to $m=9$ model parameters demonstrate speedups of $85\times$--$120\times$ relative to the full-order solver while maintaining excellent accuracy, establishing tensorial surrogates as a rigorous and efficient computational foundation for many-query workflows.
\end{abstract}

\section{Introduction}\label{s:intro}

Integrating simulation with optimization is a powerful paradigm for informed decision-making in biomedical applications~\cite{mang2020:integrated,biros2023:inverse, mang2018:pdeconst}. Even the solution of the forward problem alone is often computationally demanding. In this work, we consider a coupled, nonlinear, time-dependent parameter-to-state map whose repeated evaluation becomes prohibitive in many-query settings such as inverse problems,statistical inference, and digital twin applications~\cite{ghattas2021:learning, hogea2008:image,scheufele2020:imagedriven, ghafouri2025:inverse, willcox2021:imperative,knopoff2013:adjoint}.

To address this computational bottleneck, we investigate the construction of surrogate models based on different variants of reduced-order models ({\bf \rom}s) of the discrete forward operator~\cite{benner2015:survey, benner2017:mra, hesthaven2016:certified, rozza2008:reduced, antoulas2005:approximation}. The underlying model consists of coupled semi-linear parabolic reaction--diffusion-type equations parameterized by an $m$-dimensional parameter vector $\theta \in \mathbb{R}^m$, posed on heterogeneous two- and three-dimensional domains. The state variables represent tumor cell populations, oxygen concentration, and tissue compartments, following established tumor growth modeling frameworks~\cite{saut2014:multilayer, subramanian2019:simulation, mang2020:integrated, ghafouri2025:inverse}.

Designing effective {\rom}s that accurately capture parametric and temporal variability in large-scale, time-dependent partial differential equations (\textbf{PDE}s) remains a significant challenge, particularly in high-dimensional parameter regimes~\cite{peherstorfer2018:survey}. In the present work, we demonstrate results for problems with parameter dimensions of up to $m=9$.

\subsection{Outline of the Method}

We consider coupled, parametric, semi-linear reaction--diffusion type PDE systems of the general form
\[
\partial_t u(x,t) + \mathcal{L}(u,x,t;\theta) + \psi(u,x,t; \theta) = 0
\]

\noindent in $\Omega_B \times (0,1]$, $\Omega_B \subset \mathbb{R}^d$, $d \in \{2,3\}$, controlled by the model parameters $\theta \in \mathbb{R}^m$, with state variable $u :\Omega_B \times [0,1] \to \mathbb{R}_+^{n_s}$, initial condition $u(x,t=0) = u_0(x)$, and Neumann boundary conditions on $\partial\Omega_B$. Here, $\mathcal{L}$ is an inhomogeneous, linear diffusion operator and $\psi$ is a (typically non-linear) reaction term. We discretize this system using a finite volume method (\textbf{FVM}). For time integration, we consider first- and second-order accurate splitting schemes that treat the linear and non-linear terms individually. We consider several \rom\ variants to speed up the solution of the underlying PDE system. As a baseline, we use a proper orthogonal decomposition ROM (\textbf{\podrom}). In addition, we consider several variants of tensorial ROMs (\textbf{\trom}s)~\cite{mamonov2022:interpolatory, mueller2025:tensor, olshanskii2025:approximating, mamonov2024:tensorial, mizan2025:parametric}. The first group of {\trom}s\ uses a higher-order singular value decomposition to construct the surrogate (\textbf{\htrom})~\cite{de2000:multilinear}. The second group of {\trom}s\ uses a tensor-train decomposition (\textbf{\ttrom})~\cite{oseledets2011:tensor_train}.

\subsection{Related Work}\label{s:literature}

Works discussing tumor growth models of varying complexity are abundant in the literature, ranging from early reaction–diffusion formulations to multi-physics and multi-scale descriptions~\cite{swanson2003:virtual, swanson2011:quantifying,hawkinsdaarud2012:numerical, hormuth2017:mechanically, jarrett2018:mathematical,lima2017:selection, oden2010:general, oden2013:selection, jbabdi2005:simulation,yankeelov2013:clinically, oden2016:toward, ozuugurlu2015:note}. In the present work, we consider a parametric, coupled system of semi-linear parabolic PDEs with six state variables representing tumor cell populations, oxygen concentrations, and tissue compartments. The interactions among these state variables are governed by nine parameters $\theta$. Variants of this model have been studied previously, for example in~\cite{saut2014:multilayer, subramanian2019:simulation, mang2020:integrated, ghafouri2025:inverse}. Here, we do not account for mechanical coupling with the surrounding brain parenchyma; formulations that include such effects can be found in~\cite{hogea2008:image, gooya2012:glistr, hormuth2017:mechanically}.

These types of models have been employed in a range of application-driven settings, including image analysis workflows~\cite{gooya2012:glistr,scheufele2019:coupling, gholami2017:framework, mang2018:pdeconst,bakas2015:glistrboost, bakas2018:identifying}, treatmentplanning or design~\cite{lima2017:selection, liu2025:m4rl,sabir2017:mathematical, krawczyk2020:optimal, yadav2025:understanding}, and parameter estimation for prediction and decision support~\cite{ghafouri2025:inverse, scheufele2020:fully,scheufele2020:imagedriven, mang2012:biophys, hogea2008:image,gholami2016:inverse, zhang2025:personalized, chen2023:tgmnet,metz2024:deepgrowth, knopoff2013:adjoint}. These workflows typically require estimation of the model parameters $\theta$ from data, an aspect we do not address in the present work; instead, our focus is exclusively on the construction of efficient surrogate models for forward simulation.

Driven by substantial advances across applied mathematics, engineering, and computing hardware, there has been increasing interest in machine-learning-based surrogate models in recent years~\cite{metz2024:deepgrowth, laslo2025:diffusion, zhang2024:pinninfiltration, pati2020:estimating, weidner2025:trainforwards, zhang2025:personalized, chen2025:deep, raissi2019:physics, lu2021:learning, li2020:fourier}. We take a different, more mathematically grounded approach that provides rigorous approximation guarantees. In particular, we consider several \rom\ variants~\cite{benner2015:survey, benner2017:mra, hesthaven2016:certified, quarteroni2015:reduced}. Beyond such guarantees, projection-based {\rom}s offer further advantages over purely data-driven surrogates, including improved data efficiency, enhanced interpretability, certified error control, and predictable online computational complexity. Moreover, these approaches preserve a direct connection to the underlying governing equations, a property that is especially desirable in safety-critical and decision-support settings. A comprehensive overview of projection-based model reduction and hyper-reduction techniques can be found in~\cite{benner2017:mra}.

The use of {\trom}s in the context of parabolic equations is not new~\cite{manzini2025:low, dolgov2012:fast, mamonov2022:interpolatory, mamonov2025:priori}. Tensorial and interpolatory reduction techniques can be understood as tools for recovering efficient offline/online decompositions in the presence of nonlinearities or nonaffine parameter dependence, particularly for problems exhibiting polynomial or multilinear structure~\cite{buithanh2008:modelreduction}. In classical reduced-basis settings, such decompositions are readily available for affine operators, but become significantly more challenging once nonlinear terms are introduced. Hyper-reduction strategies, such as the empirical interpolation method and its discrete variant, address this issue by introducing additional approximation spaces for nonlinear operators~\cite{barrault2004:eim, chaturantabut2010:deim}. For problems with polynomial nonlinearities, tensorial variants of \podrom\ provide an alternative route to online efficiency by exploiting multilinear structure in the reduced operators, leading to reduced online costs that depend only on the reduced dimension~\cite{stefanescu2014:tensorialpod}. In this sense, {\trom}s offer a natural framework for capturing parametric dependence while maintaining explicit control over computational complexity as the parameter dimension $m$ increases. We note that our numerical framework avoids certain difficulties associated with nonlinearities through the use of an  operator-splitting strategy.

\subsection{Contributions}\label{s:contributions}

We design and evaluate a numerical framework for the efficient surrogate-based simulation of brain tumor growth. Our major contributions are:
\begin{itemize}[leftmargin=*]
\item
We design an effective numerical framework for solving coupled, semi-linear reaction--diffusion type PDE systems arising in brain tumor growth modeling. The framework combines a finite volume discretization with an operator-splitting time integration strategy that separates linear diffusion from nonlinear reaction terms, enabling a straightforward and efficient deployment of projection-based surrogate models.
\item
We develop and compare four variants of tensorial reduced-order models: an interpolatory \htrom\ and an interpolatory \ttrom, each available in intrusive (projection-based) and non-intrusive implementations. These surrogates exploit the multilinear structure of the snapshot tensor to achieve simultaneous compression across spatial, temporal, and parameter dimensions.
\item
We introduce a spatial pre-compression stage that projects all snapshots onto a low-dimensional spatial subspace prior to tensor decomposition, rendering the offline construction tractable for large-scale three-dimensional problems.
\item
We provide a detailed empirical evaluation of the performance of the designed surrogate models for two brain tumor growth models of varying complexity---a scalar single-species and a six-state, nine-parameter multi-species formulation---in both two and three spatial dimensions.
\item
The proposed \trom\ surrogates deliver online speedups between $85\times$ and $120\times$ relative to the full-order solver while maintaining excellent approximation accuracy, significantly outperforming the \podrom\ baseline across all test cases considered.
\end{itemize}

\subsection{Limitations}\label{s:limitations}

We identify several open issues and remaining research directions.
\begin{itemize}[leftmargin=*]
\item
Our methodology has been demonstrated for up to $m=9$ parameters; extending it to significantly higher parameter dimensions requires additional work, as the size of the training set grows exponentially with $m$ and the offline cost of constructing the snapshot tensor becomes rapidly prohibitive. To overcome this limitation, tensor completion and cross-approximation techniques have recently been studied in the context of \trom\ training~\cite{mamonov2024slice,budzinskiy2025low}; however, these approaches fall beyond the scope of the present paper.
\item
Our implementation is currently a research prototype written in \texttt{MATLAB}, which imposes practical limits on performance and scalability. Deploying the framework in a compiled, high-performance computing-ready language and exploiting parallelism in both the offline (snapshot generation) and online (tensor decomposition, reduced integration) stages is left to future work.
\item
The three-dimensional numerical experiments reported in this work are restricted to the single-species model; the multi-species formulation has only been validated in two spatial dimensions. Extending the multi-species model to three dimensions is an important next step but requires additional implementation effort and substantially larger computational resources for offline snapshot generation.
\item
The performance of the proposed surrogates in extrapolation regimes---i.e., for parameter values outside the training \textit{range} $\Theta_s$---has not been systematically studied and remains an open question. This will become of particular importance for deploying the proposed methodology in the context of inference of parameter values from clinical data~\cite{mang2020:integrated, ghafouri2025:inverse, mang2012:biophys, mang2018:pdeconst, gholami2016:inverse}.
\end{itemize}

\subsection{Outline}\label{s:outline}

We present mathematical models in \Cref{s:model}. This includes a single- and a multi-species model (see \Cref{s:single_species_model} and \Cref{s:multi_species_model}). We describe our numerical approach in \Cref{s:numerics}. We describe the discretization in \Cref{s:discretization}. We present our numerical time integration strategy in \Cref{s:time_integration}. We describe the designed surrogate models in \Cref{s:surrogate}. This includes a baseline \podrom\ (see \Cref{s:pod}), and four variants of {\trom}s---an \htrom\ and a \ttrom\ (see \Cref{s:tucker} and \Cref{s:tt}), each available in an intrusive and a non-intrusive implementation. We present numerical results in \Cref{s:results} and conclude with \Cref{s:conclusions}.

\section{Mathematical Model}\label{s:model}

We present the mathematical models we are going to consider below. The first model considers a single-species semi-linear parabolic reaction--diffusion type model for the tumor cell density. The second model consists of a coupled system of a similar type that accounts for different cell phenotypes and nutrition.

\subsection{Single Species Model}\label{s:single_species_model}

We model the tumor cell concentration $u(x,t)$ using the semi-linear parabolic PDE~\cite{mang2014:methoden, mang2012:biophys, gholami2016:inverse}
\begin{equation}\label{e:singlespecies}
\begin{aligned}
\partial_t u(x,t) - \idiv \alpha(x) \igrad u(x,t) + \rho\, u(x,t) \left( 1 - u(x,t) \right) & = 0 && \text{in}\,\,\Omega_B \times (0,1], \\
u(x,t) & = u_0(x) && \text{in}\,\,\Omega_B \times \{0\}, \\
\partial_n u(x,t) & = 0 && \text{in}\,\,\partial\Omega_B \times [0,1],
\end{aligned}
\end{equation}

\noindent where $\Omega_B \subset \mathbb{R}^d$ denotes the domain occupied by brain parenchyma with boundary $\partial\Omega_B$ and closure $\bar{\Omega}_B$. The diffusion operator in \Cref{e:singlespecies} accounts for the migration of cells into surrounding (healthy) tissue. It is controlled by the scalar map $\alpha : \bar{\Omega}_B \to \mathbb{R}_+$. The third term in the equation is a non-linear logistic reaction operator that accounts for the proliferation of cells controlled by the proliferation rate $\rho > 0$. The scalar field $\alpha$ is defined by
\[
\alpha(x)
\defeq \alpha_{\textit{w}}\pi_{\textit{w}}(x) + \alpha_{\textit{g}}\pi_{\textit{g}}(x)
= \alpha(\pi_{\textit{w}}(x) + (1/10)\pi_{\textit{g}}(x)),
\]

\noindent where $\pi_{\textit{w}} : \bar{\Omega} \to [0,1]$ and $\pi_{\textit{g}} : \bar{\Omega} \to [0,1]$ are white and gray matter probability maps obtained from a brain atlas so that $0 \le \pi_{\textit{w}}(x) + \pi_{\textit{g}}(x) \le 1$ for all $x\in\Omega_B$ and $\pi_{\textit{w}}(x) = \pi_{\textit{g}}(x) = 0$ for $x \in \Omega \setminus \Omega_B$. The migration into surrounding healthy tissue is controlled by the scalar parameters $\alpha_w = \alpha > 0$ and $\alpha_g = \alpha/10 > 0$ (we simplify the model to be controlled by a single parameter $\alpha$).

\subsection{Multi-Species Model}\label{s:multi_species_model}

We extend the single-species formulation to a mathematical framework that captures the interplay between multiple tumor cell phenotypes, oxygen dynamics, and brain tissue heterogeneity. Similar formulations have been considered in~\cite{saut2014:multilayer, subramanian2019:simulation, mang2020:integrated, ghafouri2025:inverse}. We do not include the mechanical coupling with linear elasticity equations, thereby neglecting tumor-induced mass effect.

The model is based on the \emph{go-or-grow} hypothesis, which postulates that tumor cells switch between two primary phenotypes: proliferative, denoted by $\up : \bar{\Omega}_B\times[0,1]\to\mathbb{R}_+$, and infiltrative, denoted by $\ui : \bar{\Omega}_B\times[0,1]\to\mathbb{R}_+$. In nutrient-rich and oxygenated regions, proliferative cells undergo rapid mitosis, driving local tumor expansion. Under hypoxic conditions, cells transition to the infiltrative phenotype, enabling migration toward better-oxygenated regions. Upon encountering favorable conditions, infiltrative cells revert to the proliferative phenotype, continuing the growth–migration cycle. Severe hypoxia leads to cell death and the formation of a necrotic core, denoted by $\un$. Oxygen concentration, denoted by $\uo : \bar{\Omega}_B\to\mathbb{R}_+$, plays a central regulatory role, influencing both proliferation rates and phenotypic switching. The oxygen field is consumed by tumor cells and supplied by healthy vasculature, with rates $\kappa_c \geq 0$ and $\kappa_s \geq 0$, respectively. Gray matter, denoted by $\ug : \bar{\Omega}_B\to\mathbb{R}_+$, and white matter, $\uw : \bar{\Omega}_B\to\mathbb{R}_+$, fractions contribute to the spatially varying diffusion coefficient $\alpha : \bar{\Omega} \to \mathbb{R}_+$ and proliferation rate $\rho : \bar{\Omega} \to \mathbb{R}_+$. These tissue maps remain fixed in time and satisfy $\ug + \uw = \uh$, where $\uh$ denotes the healthy tissue.

Overall, we model the evolution of the densities $u : \bar{\Omega}_B \times[0,1] \to \mathbb{R}^{n_s}_+$, with $n_s =6$,
\[
(x,t) \mapsto u(x,t) = (\up, \ui, \un, \uo, \ug, \uw)(x,t),
\]

\noindent as the following coupled, semi-linear system of parabolic PDEs defined in $\Omega_B \times (0,1]$:
\begin{subequations}\label{e:multispecies}
\begin{align}
\partial_t \up  - f(\up) + \gamma_p^u \up (1-\ui) - \gamma_i^u \ui (1-\up) + \nu_o^u \up (1-\un) &= 0,
\label{e:proliferative}
\\
\partial_t \ui  - \idiv \alpha\igrad \ui - f(\ui) + \gamma_i^u \ui (1-\up) -\gamma_p^u \up(1-\ui) + \nu_o^u \ui (1-\un) &= 0,
\label{e:invasive}
\\
\partial_t \un - \nu_o^u (\ui + \up) (1-\un) &= 0,
\label{e:necrotic}
\\
\partial_t \uo + \kappa_c \uo \up  - \kappa_s (1-\uo) (\uw+\ug) &= 0,
\label{e:oxygen}
\\
\partial_t \uw + (\uw/(\uw+\ug))\left(\idiv \alpha^u\igrad \ui + f(\up) + f(\ui)\right) & =0,
\label{e:white_matter}
\\
\partial_t \ug + (\ug/(\uw+\ug))\left(\idiv \alpha^u\igrad \ui+f(\up)+f(\ui)\right) & =0,
\label{e:gray_matter}
\end{align}
\end{subequations}

\noindent with initial conditions $ \up = \up^0 $, $ \ui^0 = 0$, $\un^0 = 0$, $\uo^0 = 1$, $\uw = \uw^0 $ and $\ug = \ug^0$ in $\Omega_B \times \{0\}$ and Neumann boundary conditions in $\partial\Omega_B \times [0,1]$. Here, the equations \Cref{e:proliferative}, \Cref{e:invasive}, \Cref{e:necrotic} model the evolution of proliferative tumor cells $\up$, infiltrative tumor cells $\ui$, and necrotic tumor cells $\un$. Equation \Cref{e:oxygen} accounts for the temporal evolution of the oxygen concentration $\uo$. The last two equations \Cref{e:white_matter} and \Cref{e:gray_matter} model the change of the healthy tissue $\ug$ (gray matter) and $\uw$ (white matter) attributed to infiltrating and proliferating tumor cells. We summarize the fields that appear in \Cref{e:multispecies} in \Cref{t:multispecies-fields}.

\begin{table}
\caption{List of variables and parameter maps used in the multispecies go-or-grow tumor model. This table summarizes the notation used for the tumor cell populations, oxygen dynamics, and transition rates between proliferative and infiltrative phenotypes. The model is based on~\cite{gholami2016:inverse,saut2014:multilayer,subramanian2019:simulation,ghafouri2025:inverse}.}
\label{t:multispecies-fields}
\tabadjust
\begin{tabular}{llc}
\toprule
\textbf{symbol} & \textbf{description}  & \textbf{equation}\\
\midrule
$\up(x, t)$ & density of proliferative tumor cells & \Cref{e:proliferative}\\
$\ui(x, t)$ & density of infiltrative tumor cells & \Cref{e:invasive} \\
$\un(x, t)$ & density of necrotic tumor cells & \Cref{e:necrotic}\\
$\uo(x, t)$ & oxygen concentration in the tumor micro-environment & \Cref{e:oxygen} \\
$\ug(x, t)$ & density of gray matter cells & \Cref{e:gray_matter} \\
$\uw(x, t)$ & density of white matter cells & \Cref{e:white_matter} \\
\midrule
$\alpha(x,t)$ & diffusion coefficient map & \Cref{e:diffusion-map} \\
$f(u)$ & net reaction term modeling tumor proliferation  & \Cref{e:net-reaction} \\
$\rho^u(x,t)$ & proliferation rate &  \Cref{e:prolif-rate}\\
$\gamma_p^u(x, t)$ & transition rate for proliferative to infiltrative phenotype & \Cref{e:transition-pro} \\
$\gamma_i^u(x, t)$ & transition rate for infiltrative to proliferative phenotype & \Cref{e:transition-inv} \\
$\nu_o^u(x, t)$ & oxygen-dependent transition threshold & \Cref{e:oxygen-death} \\
\bottomrule
\end{tabular}
\end{table}

The model in \Cref{e:multispecies} is controlled by various parameters. We summarize these parameters in \Cref{t:multispecies-para}. We introduce a superscript $u$ if the value of the parameter depends on the species that appear in our model. In the following, we derive principled rules for computing these parameters and introduce suitable simplifications to reduce the size of the parameter space, with a view toward its use in an inverse problem.

The field $\alpha^u : \bar{\Omega}_B \to \mathbb{R}_{+}$ that appears in~\Cref{e:white_matter} and~\Cref{e:gray_matter} controls the infiltration of $\ui$ into surrounding healthy tissue. We assume that infiltration is different in white matter than in gray matter. The model for $\alpha$ is given by
\begin{equation}\label{e:diffusion-map}
\alpha^u(x,t)
= \alpha_w \uw(x,t) + \alpha_g \ug(x,t)
= \alpha (\uw(x,t) + (1/5)\ug(x,t))
\end{equation}

\noindent with parameters $\alpha_w = \alpha > 0$ and $\alpha_g = (\alpha/5) > 0$ (we reduce the model to be controlled by a single parameter $\alpha > 0$).

The net reaction term $f$ that appears in \Cref{e:proliferative}, \Cref{e:invasive}, \Cref{e:oxygen}, respectively, accounts for logistic proliferation modulated by oxygen availability. The model is controlled by the proliferation rate $\rho^u : \bar{\Omega}_B \times [0,1] \to \mathbb{R}_+$, and the thresholds $\uo^{\textit{inv}} > 0$, $\uo^{\textit{mit}} > 0$, and $\uo^{\textit{hyp}} > 0$ for invasion, mitosis, and hypoxia, where $\uo^{\textit{mit}} = (\uo^{\textit{hyp}} + \uo^{\text{inv}})/2$. Let $u_{\textit{s}}(x,t) = \up(x,t) + \ui(x,t) + \un(x,t)$ for all $t\in[0,1]$, $x\in\Omega_B$, with $0 \le u_{\textit{s}} \le 1$. Then, for some concentration $\tilde{u}(x,t)$, we have
\begin{equation}\label{e:net-reaction}
f(\tilde{u})(x,t) =
\begin{cases}
\rho^{\tilde{u}}(x,t) \tilde{u}(x,t)(1 - u_{\textit{s}}(x,t)) & \text{for}\,\, \uo > \uo^{\textit{inv}}, \\
\rho^{\tilde{u}}(x,t) \tilde{u}(x,t)(1 - u_{\textit{s}}(x,t)) \dfrac{\uo - \uo^{\textit{mit}}}{\uo^{\textit{inv}} - \uo^{\textit{mit}}} &\text{for}\,\, \uo^{\textit{inv}} \geq \uo \geq \uo^{\textit{mit}}, \\
0 & \text{for}\,\, \uo < \uo^{\textit{mit}}.
\end{cases}
\end{equation}

The factor $(1 - u_{\textit{s}}(x,t))$ ensures that there will be no reaction or increase of tumor concentration once the total tumor concentration $u_{\textit{s}} \approx 1$. We model $\rho^u$ as
\begin{equation}\label{e:prolif-rate}
\rho^{\tilde{u}}(x, t) = \rho_w^{\tilde{u}} \uw(x, t) + \rho_g^{\tilde{u}}\ug(x, t)
\end{equation}

\noindent controlled by the parameters $\rho_g^{\tilde{u}} > 0$ and $\rho_w^{\tilde{u}} > 0$. These parameters depend on the species $u$. To reduce the dimension of our parameter space, we collapse the model to be controlled by a single parameter $\rho > 0$. For $\up$ we set
\[
\rho^{\up}(x, t) = \rho (\uw(x, t) + (1/5) \ug(x, t)),
\]

\noindent i.e., $\rho_w^{\up} = \rho$ and $\rho_g^{\up} = \rho/5$. Likewise, for $\ui$ we choose
\[
\rho^{\ui}(x, t) = (\rho/10) (\uw(x, t) + (1/5) \ug(x, t)),
\]

\noindent i.e., $\rho_w^{\ui} = \rho/10$ and $\rho_g^{\ui} = \rho/50$.

Oxygen-dependent cell death is controlled by the parameter $\nu_o^u$. This parameter appears in \Cref{e:proliferative}, \Cref{e:invasive}, \Cref{e:necrotic}. We model it as
\begin{equation}\label{e:oxygen-death}
\nu_o^u(x,t) = \lambda_d h_\varepsilon(\uo^{\textit{hyp}} - \uo(x,t)),
\end{equation}

\noindent where $h_{\varepsilon}(x) = 1/(1 + \exp(-x/\varepsilon))$ is a smooth approximation of the Heaviside function, where $\varepsilon > 0$ controls the smoothness of the transition and $\lambda_d > 0$ is the death rate of the tumor cells.

The switching between phenotypes $\ui$ and $\up$ is controlled by the transition parameters $\gamma_p^u$ and $\gamma_i^u$ that appear in \Cref{e:proliferative} and \Cref{e:invasive}. We model the rates according to
\begin{subequations}
\begin{align}
\gamma_p^u &= \gamma_p h_{\varepsilon}(\uo^{\textit{inv}} - \uo), \label{e:transition-pro} \\
\gamma_i^u &= \gamma_i h_{\varepsilon}(\uo - \uo^{\textit{inv}}), \label{e:transition-inv}
\end{align}
\end{subequations}

\noindent controlled by the scalar parameters $\gamma_p > 0$ and $\gamma_i > 0$. The remaining user defined parameters are $\kappa_s > 0$ (the oxygen supply rate) and $\kappa_c > 0$ (the oxygen consumption rate); they appear in~\Cref{e:oxygen}. In summary, based on the simplifications that allowed us to reduce the dimension of the parameter space $\Theta \subset \mathbb{R}^m_{+}$, the considered model is controlled by the following $m=9$ parameters
\[
\theta = (\alpha, \rho,\gamma_p, \gamma_i, \lambda_d, \kappa_s, \kappa_c, \uo^\textit{inv}, \uo^\textit{hyp}) \in \Theta.
\]

\noindent We summarize these parameters, their meaning, and values and parameter ranges considered in the literature in \Cref{t:multispecies-para}.

\begin{table}
\caption{Model parameters, values, and parameter ranges for the system \Cref{e:multispecies}. The parameters correspond to reaction and diffusion coefficients, go-or-grow transition rates, oxygen dynamics, and screening properties of the model considered and developed in~\cite{saut2014:multilayer,subramanian2019:simulation,ghafouri2025:inverse}. The values included here are based on~\cite{gholami2016:inverse,saut2014:multilayer}. Derived values are highlighted in gray.}
\label{t:multispecies-para}
\tabadjust
\begin{tabular}{llccc}\toprule
\textbf{symbol} & \textbf{description} & \textbf{default value} & \textbf{range} & \textbf{equation}\\
\midrule
$\alpha$                             & diffusion coefficient & $0.05$ & [\fnum{1e-3}, \fnum{1e-1}]  & \Cref{e:diffusion-map}\\
\rowcolor{gray!40}$\alpha_w$         & diffusion coefficient in white matter & $\alpha$ &   & \\
\rowcolor{gray!40}$\alpha_g$         & diffusion coefficient in gray matter  & $0.2\alpha = 0.01$ &  & \\
$\rho$                               & reaction coefficient & $18.0$ & [\fnum{5.0}, \fnum{2.5e1}] & \Cref{e:prolif-rate} \\
\rowcolor{gray!40}$\rho_w^{\up}$     & reaction coefficient in white matter & $\rho$  &  & \\
\rowcolor{gray!40}$\rho_g^{\up}$     & reaction coefficient in gray matter & $0.2\rho = 3.6$ &  & \\
\rowcolor{gray!40}$\rho_w^{\ui}$     & reaction coefficient in white matter & $0.1\rho = 1.8$ &   & \\
\rowcolor{gray!40}$\rho_g^{\ui}$     & reaction coefficient in gray matter & $0.02\rho = 0.36$ &   & \\
$\gamma_{p}$                         & transition rate from $\up$ to $\ui$ & $5.0$ & [\fnum{1e-1}, \fnum{1.0e1}] & \Cref{e:transition-pro}\\
$\gamma_{i}$                         & transition rate from $\ui$ to $\up$ & $7.0$ & [\fnum{1e-1}, \fnum{1.5e1}] & \Cref{e:transition-inv}\\
$\lambda_d$                          & death rate & $10.0$ & [\fnum{1.0}, \fnum{2.0e1}] & \Cref{e:oxygen-death}\\
$\uo^{\textit{hyp}}$                 & hypoxia oxygen threshold & $0.5$ & [\fnum{1e-3}, \fnum{6e-1}] & \Cref{e:oxygen-death} \\
$\uo^{\textit{inv}}$                   & invasive oxygen threshold  & $0.90$ & [\fnum{7e-1}, \fnum{1.0}] & \Cref{e:transition-pro},\Cref{e:transition-inv}\\
$\kappa_s$ & oxygen supply rate      & $6$  & [\fnum{1.0}, \fnum{8.0}] & \Cref{e:oxygen}\\
$\kappa_c$ & oxygen consumption rate & $15.0$ & [\fnum{1.0}, \fnum{2.0e1}] & \Cref{e:oxygen}\\
\bottomrule
\end{tabular}
\end{table}

\section{Numerical Methods}\label{s:numerics}

In this section, we introduce a numerical approach for solving the dynamical systems \Cref{e:singlespecies} and \Cref{e:multispecies}. For simplicity of presentation, we develop the methodology using the prototypical equation,
\begin{equation}\label{e:prototype-numerics}
\partial_t u(x,t) - \idiv \alpha(x) \igrad u(x,t) + \psi(u,x,t) = 0
\quad\text{in}\quad\Omega_B \times (0,1]
\end{equation}

\noindent with initial condition $u(x,t) = u_0(x)$ in $\Omega_B \times \{0\}$ and Neumann's boundary conditions $\partial_n u(x,t) = 0$ in $\partial\Omega_B \times [0,1]$. This equation allows us to consolidate \Cref{e:singlespecies} and \Cref{e:multispecies} into a single model; the relation to \Cref{e:singlespecies} is straightforward; for \Cref{e:multispecies}, we can assume with a slight abuse of notation that $u(x,t) = (\up, \ui, \un, \uo, \ug, \uw)(x,t)$. To keep the notation light, we mostly limit the description of the numerical approach to a single species model. For the ROMs considered in our work, we construct reduced models for each species as outlined below.

\subsection{Numerical Discretization}\label{s:discretization}

We discretize the equations using a finite volume method~\cite{eymard2000finite, leveque2002:fvm}. Our discretization strategy follows~\cite{mang2014:methoden}. Let $d\in\{1,2,3\}$ denote the dimension of the ambient space. We discretize the function $u : \bar{\Omega} \to \mathbb{R}$ with $\Omega \defeq [\omega_1,\omega_2] \times \cdots \times [\omega_{2d-1},\omega_{2d}] \subset \mathbb{R}^d$ on a cell centered grid $\Omega^h\in \mathbb{R}^{d \times n_1 \times \cdots \times n_d}$, where $n_i \in \mathbb{N}$ denotes the number of cells. The grid coordinates $x_k = (x_k^1,\ldots,x_k^d) \in \mathbb{R}^d$, $k \in \mathbb{Z}^d$, with $x_k^i = \omega_{2i-1} + (k_i - 1/2)h_i$ represent the center point of a cell of width $h = (h_1,\ldots,h_d) \in \mathbb{R}^d$, $h_i = (\omega_{2i} - \omega_{2i-1})/n_i$. Similarly, we discretize the temporal domain $[0,1]$ using a nodal mesh consisting of $n_t + 1$ time points with mesh size $h_t = 1/n_t$ and time points $t^j = (j-1)h_t$. Consequently, the numerical approximation of $u$ at location $x_k$ and time point $t^j$ is given by $u_k^j \approx u(x_k, t^j)$. We arrange the discrete $u_k^j$ in lexicographical ordering to obtain
\[
u^j_h = (u_{1,\ldots,1}^j,\ldots, u_{n_1,\ldots,n_d}^j) \in \mathbb{R}^{n_x},
\]

\noindent where $n_x = \prod_{i=1}^d n_i$ denotes the total number of mesh points.

Before we introduce our scheme for integrating the equation in time, we are going to discuss the discretization in space. As such, we will work with the semi-discrete state variable
\[
u_h(t) = (u_{1,\ldots,1}(t),\ldots, u_{n_1,\ldots,n_d}(t)) \in \mathbb{R}^{n_x}.
\]

We use short differences to discretize the gradient operator. We obtain the sparse matrix
\[
G_i = \frac{1}{h^i}
\begin{pmatrix}
b_1 & b_2  &        &       &     &     \\
    & -1   & 1      &       &     &     \\
    &      & \ddots & \ddots&     &     \\
    &      &        & -1    &   1 &     \\
    &      &        &       & b_3 & b_4 \\
\end{pmatrix} \in \mathbb{R}^{n_i + 1 \times n_i}
\]

\noindent for the $i$th coordinate direction; $b_j$ represent the entries for the boundary conditions. The two and three-dimensional gradient operators are constructed using Kronecker products; we have
\[
\operatorname{GRAD} =
\begin{pmatrix}
 I_{n_1} \otimes G_1\\
 G_2 \otimes I_{n_2}
\end{pmatrix} \in \mathbb{R}^{(n_1 + 1)n_2 + (n_2 + 1)n_1 \times n_1n_2}
\]

\noindent for $d=2$ and
\[
\operatorname{GRAD} =
\begin{pmatrix}
I_{n_3} \otimes I_{n_2} \otimes G_1 \\
I_{n_3} \otimes G_2 \otimes I_{n_1} \\
G_3 \otimes I_{n_2} \otimes I_{n_1} \\
\end{pmatrix} \in \mathbb{R}^{\tilde{n} \times n_1n_2n_3}
\]

\noindent with $\tilde{n} = (n_1 + 1)n_2n_3 + n_1(n_2 + 1)n_3 + n_1n_2(n_3 + 1)$. Notice that for $d=1$, $\operatorname{GRAD}$  maps from a cell-centered grid to a nodal grid. For $d=2$ and $d=3$, $\operatorname{GRAD}$ maps from a cell-centered grid to a staggered grid.

The divergence operator is given by
\[
\operatorname{DIV} =
\begin{pmatrix}
 I_{n_1} \otimes G_1 &  G_2 \otimes I_{n_2}
\end{pmatrix} \in \mathbb{R}^{n_1n_2 \times (n_1 + 1)n_2 + (n_2 + 1)n_1}
\]

\noindent for $d=2$ and
\[
\operatorname{DIV} =
\begin{pmatrix}
I_{n_3} \otimes I_{n_2} \otimes G_1 &
I_{n_3} \otimes G_2 \otimes I_{n_1} &
G_3 \otimes I_{n_2} \otimes I_{n_1}
\end{pmatrix} \in \mathbb{R}^{n_1n_2n_3 \times \tilde{n}}
\]

\noindent for $d=3$. The operator $\operatorname{DIV}$ maps back to the original cell-centered grid.

We assume that the coefficient map $\alpha$ is also allocated on a cell-centered grid. To map $\alpha_h \in \mathbb{R}^{n_x}$ to a staggered grid, we introduce the one-dimensional grid change operator
\[
M_i = \frac{1}{2}
\begin{pmatrix}
3 & -1     &        &   \\
1 &  1     &        &   \\
  & \ddots & \ddots &   \\
  &        &   1    & 1 \\
  &        &  -1    & 3 \\
\end{pmatrix} \in \mathbb{R}^{n_i + 1 \times n_i}.
\]

\noindent The $d$-dimensional operator $\operatorname{AVG}$ is constructed with the same stencil as $\operatorname{GRAD}$. Using these matrices, we can approximate the diffusion operator $\mathcal{L}(\alpha) = -\idiv \alpha(x) \igrad$ as
\[
L_h \defeq -\operatorname{DIV} \operatorname{diag}(e \oslash (\operatorname{AVG} \alpha_h)) \operatorname{GRAD} \in \mathbb{R}^{n_x \times n_x},
\]

\noindent where $e = \begin{pmatrix} 1,\ldots,1\end{pmatrix} \in \mathbb{R}^{\tilde{n}}$, $\oslash : \mathbb{R}^{\tilde{n}} \to \mathbb{R}^{\tilde{n}}$ denotes element-wise division (Hadamard division) and $e \oslash (\operatorname{AVG} \alpha_h)$ implements harmonic averaging.

The discretization of $\psi$ is straightforward; for example, for $\psi(u,x,t) =  \rho u(x,t) (1 - u(x,t))$, we obtain the semi-discrete expression
\[
\psi_h(u_h(t),t) = \rho u_h(t) \odot (e - u_h(t)) \in \mathbb{R}^{n_x}.
\]

\noindent Overall, we arrive at the semi-discrete system of ODEs
\begin{equation}\label{e:semi-discretized-ode}
\text{d}_t u_h(t) = f_h(u_h(t),t)
\end{equation}

\noindent with $f_h(u_h(t),t) \defeq L_h u_h(t) - \psi_h(u_h(t),t)$.

One issue we have not addressed so far is that we actually need to implement these operators on a non-rectangular, complex domain $\Omega_B$ with a curved boundary $\partial\Omega_B$. For this purpose, we use a fictitious domain method to implement the boundary conditions~\cite{mang2012:biophys, mang2014:methoden, hogea2008:image}.

\subsection{Numerical Time Integration}\label{s:time_integration}

Our approach for integrating \Cref{e:prototype-numerics} in time follows ideas presented in~\cite{saut2014:multilayer, subramanian2019:simulation, gholami2016:inverse}. We consider a second order accurate Strang operator splitting scheme in the present work (cf., e.g.,~\cite{strang1968:difference,hundsdorfer2003:numerical}), splitting the equation into diffusion and reaction steps; see \Cref{a:strang-splitting}.%

\begin{algorithm}
\caption{Second order accurate Strang splitting. We solve the diffusion step using a Crank--Nicolson method. The reaction step is solved analytically or numerically (using a second-order Runge--Kutta method).\label{a:strang-splitting}}
\begin{algorithmic}[1]
\State $\hat{u}_h \gets$ solve $(I_{n_x} - (h_t/4)L)\hat{u}_h =  (I_{n_x} + (h_t/4)L)u^j_h$
\State $\tilde{u}_h \gets$ solve $\partial_t u = -\psi(u,x,t)$ in $[t^j,t^{j+1}]$ with initial condition $\hat{u}_h$
\State $u_h^{j+1} \gets$ solve $(I_{n_x} - (h_t/4)L)u_h^{j+1} =  (I_{n_x} + (h_t/4)L)\tilde{u}_h$
\end{algorithmic}
\end{algorithm}

For the reaction step (step 2), we either use an analytic expression to evaluate the time integration (when possible) or a second-order accurate Runge--Kutta scheme. We solve the linear system using a preconditioned conjugate gradient method with a tolerance of 1e--8.\footnote{We tested a first order accurate Lie splitting strategy and explored the use of BDF1 and BDF2 for solving the linear system as alternatives. Balancing numerical accuracy and computational complexity, we settled for the proposed scheme.} We precondition the Crank--Nicolson step using an incomplete Cholesky factorization with a drop tolerance of 1e--3 and a diagonal compensation with a perturbation of 1e--3. We also explored the use of a full Cholesky factorization but this strategy was not viable in the three-dimensional case because of memory restrictions. Designing a preconditioner that is faster to assemble and apply requires more work.

\subsection{Surrogate Models}\label{s:surrogate}

In general, we have to design the surrogate models for a vector valued state variable $u : \bar{\Omega}_B \times[0,1] \to \mathbb{R}^{n_s}_+$. To simplify notation, we assume that $n_s = 1$ below. Notice that we build surrogates for each individual component of $u$ independently. Consequently, the methodology described below generalizes in a straightforward way.

We denote by $u_k^j(\theta) \approx u(x_k,t^j;\theta)$ a discrete approximation of the state variable $u : \bar{\Omega}_B \times[0,1] \to \mathbb{R}_+$ at time $t^j$ and location $x_k$ as a function of the model parameters
\[
\theta = (\theta_1,\ldots,\theta_m)\in \Theta,
\quad \Theta = \bigotimes_{i=1}^m[\theta_i^{\text{min}},\theta_i^{\text{max}}] \subset \mathbb{R}^m_+.
\]

For each surrogate model considered below, we generate a training dataset $\Theta_s$ given by a Cartesian grid with cardinality $n_{\theta} = \prod_{i=1}^m n_{\theta_i}$, by sampling $n_{\theta_i} \in \mathbb{N}$ values for each parameter $\theta_i$ that controls the solution of the dynamical system within a given interval $\Theta_i \defeq [\theta_i^{\text{min}},\theta_i^{\text{max}}] \subset\mathbb{R}$ of admissible parameter values. We define the multi-index set
\begin{equation}\label{e:multiindexset-para}
\mathcal{I}_s = \left\{l = (l_1,\ldots,l_{m})  : l_i = 1,\ldots,n_{\theta_i}, \,\,i = 1,\ldots,m \right\}.
\end{equation}

The resulting dataset $\Theta_s \subset \Theta \subset \mathbb{R}^m$ of trial model parameters $\theta^l$ is given by
\[
\Theta_s = \left\{ \theta^l = (\theta_1^l, \ldots, \theta_m^l) :  \theta_i^l \in  \Theta_i,\,\, i = 1,\ldots,m,\,\, l \in \mathcal{I}_s \right\}.
\]

Based on these samples, we compute snapshots
\[
u^j_h(\theta^l) \in \mathbb{R}^{n_x},\quad
\theta^l \in \Theta_s, \quad
j = 1,\ldots,n_t,\,\, l \in \mathcal{I}_s,
\]

\noindent where $n_x = \prod_{i=1}^d n_i$ and $u^j_h(\theta^l)$ denotes a high-fidelity solution of the \fom\ (either \Cref{e:singlespecies} or \Cref{e:multispecies}) for a parameter sample $\theta^l \in \Theta_s$.

\subsubsection{Rank Selection}\label{s:rank-selection}

For the construction of the surrogate models, we compute rank-$r$ approximations. To determine the target rank $r \leq k$, $k \defeq \min\{\tilde{n},\tilde{m}\}$ of a given $\tilde{n} \times \tilde{m}$ matrix, we consider the proportion of energy $\eta(r)$ captured by the rank $r$ approximation. Let $\{\sigma_i\}_{i=1}^k$ denote the singular values of the considered matrix with $\sigma_1 \ge \sigma_2 \ge \ldots \ge \sigma_k$. Then,
\begin{equation}\label{e:energy-bound-r}
r(\epsilon) = \min\left\{ q \geq 1 : \eta(q) = \frac{\sum_{i=1}^q\sigma_i^2}{\sum_{i=1}^{k} \sigma_i^2} \geq 1-\epsilon, \,\,  k \defeq \min\{\tilde{n},\tilde{m}\}\right\},
\end{equation}

\noindent with user-defined threshold $\epsilon \in (0,1)$. For efficiency, we compute the truncated (thin) singular value decomposition (\textbf{SVD}) required in our numerical scheme using a randomized SVD (\textbf{rSVD})~\cite{halko2011:finding}. We provide a pseudo-code for the rSVD algorithm used to compute the truncated SVDs for constructing reduced basis vectors in \Cref{a:rsvd-fixed-accuracy}.

\subsubsection{Proper Orthogonal Decomposition ROM}\label{s:pod}

POD is a classical model reduction technique used to extract the most energetic modes from a set of data. It provides an optimal orthonormal basis in the least-squares sense that captures the dominant behavior of the system with significantly reduced computational cost~\cite{kunisch2002:galerkin, benner2015:survey, ghattas2021:learning}. Our goal is to find an orthonormal basis $\{\phi_i(x_k)\}_{i = 1}^r$ with truncation rank $r$ for some ansatz function such that
\[
u_k^j(\theta) \approx \sum_{i = 1}^r \mu_i(t^j,\theta) \phi_i(x_k).
\]

\paragraph{Offline Stage}
We build the POD basis during the so-called \emph{offline stage}. We first construct a snapshot matrix
\[
S =
\begin{pmatrix}
  u^1_h(\theta^{1,\ldots,1})
& u^2_h(\theta^{1,\ldots,1})
& \ldots
& u^{n_t}_h(\theta^{n_{\theta_1},\ldots,n_{\theta_m}})
\end{pmatrix}
\in \mathbb{R}^{n_x \times n_t n_{\theta}}.
\]

\noindent Subsequently, we use eigenvalues and eigenvectors of the associated Gramian matrix to find a low-dimensional POD basis $\{\phi_i\}_{i=1}^r$ that solves
\[
\begin{aligned}
\minopt_{\{\phi_i\}_{i=1}^r} &\quad \sum_{l \in \mathcal{I}_s} \sum_{j=1}^{n_t} \|u^j_h(\theta^l) -  \sum_{i = 1}^r \langle \phi_i, u^j_h(\theta^l)\rangle \phi_i \|^2_2\\
\text{subject to} & \quad \langle \phi_i, \phi_j \rangle = \delta_{ij}.
\end{aligned}
\]

\noindent We can rewrite the objective function as $\|S - \Phi\Phi^\mathsf{T}S\|_F^2$ with $\Phi = \begin{pmatrix}\phi_1 & \ldots & \phi_r \end{pmatrix} \in \mathbb{R}^{n_x \times r}$, $\Phi^\mathsf{T}\Phi = I_r$. Solving this optimization problem for $\Phi$ is equivalent to computing the \emph{compact} SVD $S = U_r \Sigma_r V_r^\mathsf{T}$, where $\Phi = U_r \in \mathbb{R}^{n_x \times r}$ denotes the POD basis (left singular vectors), $\Sigma = \operatorname{diag}(\sigma_1,\ldots,\sigma_r) \in \mathbb{R}^{r \times r}$ are the singular values, and $V_r\in\mathbb{R}^{(n_t+1)n_{\theta} \times r}$ are the right singular vectors. In practice, we do not know the value for the target rank $r\ll \min\{n_x,(n_t+1)n_{\theta}\}$ a priori; we determine $r$ based on \Cref{e:energy-bound-r} with tolerance $\epsilon^{\textit{pod}}>0$.

\paragraph{Online Stage}
For the \emph{online stage}, we use the reduced basis $\Phi = \begin{pmatrix}\phi_1 & \ldots & \phi_r \end{pmatrix} \in \mathbb{R}^{n_x \times r}$ to solve the associated reduced problem efficiently. Instead of solving for the full state $u^j_h(\theta_i) \in \mathbb{R}^{n_x}$, we solve the underlying equation in the reduced space. We consider the semi-discrete system of ODEs in \Cref{e:semi-discretized-ode}. Projecting this equation to the reduced space yields $\Phi^\mathsf{T} \text{d}_t u_h(t; \theta) = \Phi^\mathsf{T}f_h(u_h(t; \theta), \theta,x_k,t)$. Using the fact that $u_h(t;\theta) \approx \Phi \mu_h(t; \theta)$ with unknown reduced coefficients $\mu_h(t; \theta) \in \mathbb{R}^r$, we obtain the reduced dynamical system
\begin{equation}\label{e:reduced-dynamicalsys}
\text{d}_t \mu_h(t; \theta)
= f_h(\mu_h(t; \theta),\theta,x_k,t)
\end{equation}

\noindent with initial condition $\mu_h(0; \theta) = \Phi^{\mathsf{T}}u_h(0; \theta)$. Solving \Cref{e:reduced-dynamicalsys} forward in time (using our preferred method for numerical time integration) yields $\mu_h$ for all $t$. To reconstruct the full state $u_h$ at time $t$, we evaluate $u_h(t;\theta) \approx \Phi \mu_h(t;\theta)$ for all $t$.

Notice that we evaluate the non-linear terms that appear in $f$ in \Cref{e:reduced-dynamicalsys} by projecting back to the full-order space, i.e., we apply $\Phi$ to $\mu_h(t;\theta)$. Since we use an operator splitting strategy (see \Cref{s:time_integration}) and the non-linear terms are quick to evaluate, the additional computational costs associated with this step are small. As an alternative, one can consider techniques such as the \emph{discrete empirical interpolation method}~\cite{chaturantabut2010:nonlinear}. We did not explore this approach.

Next, we consider several variants of tensorial {\rom}s ({\bf \trom}). We refer to~\cite{mamonov2022:interpolatory, mueller2025:tensor, olshanskii2025:approximating, mamonov2024:tensorial, mizan2025:parametric, kilmer2021:tensor} for additional insights.

\subsubsection{Interpolatory Tensor Tucker Decomposition}\label{s:tucker}

The family of Tensor Tucker Decompositions~\cite{tucker1966:some}, in particular, the Higher-Order Singular Value Decomposition (\textbf{HOSVD})~\cite{de2000:multilinear}, is a generalization of the matrix SVD to higher-order tensors. It provides a low-rank approximation of multidimensional data, making it suitable for reducing large snapshot tensors generated by parametric dynamical systems. HOSVD offers several advantages over traditional POD by preserving the multilinear structure of snapshot data and allowing dimensionality reduction across multiple modes~\cite{de2000:multilinear, mamonov2022:interpolatory}.

\paragraph{Offline Stage}
We start by describing the \emph{offline stage} of the \htrom. For this approach, we represent the snapshot matrix $S$ as a tensor $\mathcal{S} \in \mathbb{R}^{n_x \times n_t \times n_{\theta_1} \times \cdots \times n_{\theta_m}}$ of order $m + 2$.\footnote{As stated above, our state vector has up to $n_s = 6$ components; we treat each component as a separate tensor. Again, we omit this detail from the description of the methodology to simplify the notation. Alternatively, one could assemble the components in lexicographical ordering (i.e., a ``long vector''). We tested this strategy and observed a reduction in reconstruction accuracy.} Let $\times_k$ denote the mode-$k$ tensor-matrix multiplication~\cite{de2000:multilinear}; that is, $Z = X \times_k Y$ for
\[
X = (x_{i_1i_2\ldots i_l})_{i_1,i_2,\ldots,i_l=1}^{m_1,m_2,\ldots,m_l} \in \mathbb{R}^{m_1 \times m_2 \times \cdots \times m_l}
\quad\text{and}\quad
Y = (y_{ij})_{i,j=1}^{r,m_k} \in \mathbb{R}^{r \times m_k}
\]

\noindent yields a tensor $Z$ of size $m_1 \times \cdots \times m_{k-1} \times r \times m_k \times \cdots \times m_l$ with entries
\[
z_{i_1 \ldots i_{k-1} j i_{k+1} \ldots i_l} = \sum_{i_k=1}^{m_k} x_{i_1i_2\ldots i_k\ldots i_l} y_{ji_k}.
\]

\noindent Our goal is to find orthonormal matrices $\Phi^{(x)}$, $\Phi^{(t)}$, $\{\Phi^{(\theta_i)}\}_{i=1}^m$, and a core tensor $\mathcal{C} \in \mathbb{R}^{r_x\times r_t \times r_{\theta_1} \times \cdots \times r_{\theta_m}}$ such that
\begin{equation}\label{e:hosvd-lr-snapshot}
\mathcal{S} \approx \tilde{\mathcal{S}}=
\mathcal{C} \times_1 \Phi^{(x)} \times_2 \Phi^{(t)}
\times_3 \Phi^{(\theta_1)} \ldots \times_{m+2} \Phi^{(\theta_m)}.
\end{equation}

\noindent The HOSVD delivers an efficient compression of $\mathcal{S}$ if the size of the core tensor $\mathcal{C}$ is much smaller than $\mathcal{S}$. The matrices
\[
      \Phi^{(x)} = (\phi^{(x)}_{ij})_{i,j=1}^{n_x,r_x} \in \mathbb{R}^{n_x \times r_x},
\qquad \Phi^{(t)} = (\phi^{(t)}_{ij})_{i,j=1}^{n_t,r_t} \in \mathbb{R}^{n_t \times r_t},
\qquad \Phi^{(\theta_j)} = (\phi^{(\theta_j)}_{ij})_{i,j=1}^{n_{\theta_j},r_{\theta_j}} \in \mathbb{R}^{n_{\theta_j} \times r_{\theta_j}},
\]

\noindent of the low-rank HOSVD approximation in \Cref{e:hosvd-lr-snapshot} denote the spatial modes, the temporal modes, and the modes for the $j$th entry $\theta_j \in \mathbb{R}$ of the parameter vector $\theta \in \mathbb{R}^{m}$, respectively; $r_x$, $r_t$ and $r_{{\theta}_j}$, $j=1,\ldots,m$, are the associated truncation ranks (Tucker ranks). Using \Cref{e:hosvd-lr-snapshot}, we can approximate the solution $u_k^j(\theta^l)$ for any $\theta^l = (\theta_1^l,\ldots,\theta_m^l) \in \Theta_s$ as
\begin{equation}\label{e:hosvd-stateapprox}
u_k^j(\theta^l) \approx
\sum_{i_1=1}^{r_x} \sum_{i_2=1}^{r_t}
\sum_{i_3=1}^{r_{\theta_1}} \cdots \sum_{i_{m+2}=1}^{r_{\theta_m}}
\mathcal{C}_{i_1i_2\ldots i_{m+2}}
\phi^{(x)}_{\pi(k)i_1}\,
\phi^{(t)}_{ji_2}\,
\prod_{q=1}^m \phi^{(\theta_q)}_{l_qi_{q+2}},
\end{equation}

\noindent where $\pi : \mathbb{N}^d \to \mathbb{N}$ denotes a mapping of the multi-index $k$ into lexicographical ordering.

To find the low-rank HOSVD approximation in \Cref{e:hosvd-lr-snapshot} we follow the standard algorithm introduced in~\cite{de2000:multilinear}. We compute truncated SVDs of the unfolded snapshot tensor $\mathcal{S}$ for each mode. In particular, we unfold the tensor $\mathcal{S}$ to obtain the matrices
\[
S_{(1)} \in \mathbb{R}^{n_x \times n_t n_{\theta_1} \cdots\, n_{\theta_m}}, \,
S_{(2)} \in \mathbb{R}^{n_t \times n_x n_{\theta_1} \cdots\, n_{\theta_m}}, \,
\ldots, \,
S_{(m+2)} \in \mathbb{R}^{n_{\theta_m} \times n_x n_t n_{\theta_1} \cdots\, n_{\theta_{m-1}}}.
\]

\noindent Each mode-$i$ unfolding $S_{(i)}$ reorders the entries so that each mode-$i$ fiber becomes a column of $S_{(i)}$. Subsequently, we compute the truncated (thin) SVD to obtain
\[
S_{(i)} = U^{(i)}_{r_i} \Sigma^{(i)}_{r_i}(V^{(i)}_{r_i})^{\mathsf{T}}, \quad i = 1,\dots,m+2,
\]

\noindent and set $\Phi^{(x)} = U^{(1)}_{r_1}$, $\Phi^{(t)} = U^{(2)}_{r_2}$, and $\Phi^{(\theta_j)} = U^{(j+2)}_{r_{j+2}}$, $j=1,\ldots,m$. To compute the core tensor $\mathcal{C}$, we project $\mathcal{S}$ onto the mode subspaces; that is,
\begin{equation}\label{e:core-tensor}
\mathcal{C} = \mathcal{S} \times_1 (\Phi^{(x)})^{\mathsf{T}} \times_2 (\Phi^{(t)})^{\mathsf{T}}
\times_3 (\Phi^{(\theta_1)})^{\mathsf{T}} \cdots \times_{m+2} (\Phi^{(\theta_m)})^{\mathsf{T}}
\in \mathbb{R}^{r_x \times r_t \times r_{\theta_1} \times \cdots \times r_{\theta_m}}.
\end{equation}

Each mode corresponds to a physical direction of variation (space, time, and parameters), and the core tensor $\mathcal{C}$ encodes interactions between these reduced dimensions. We note that unfolding $\mathcal{S}$ along the first dimension and computing a truncated SVD of the unfolded matrix corresponds to the POD model described in \Cref{s:pod}.

\paragraph{Online Stage} The \emph{online stage} of the \htrom\ aims to rapidly approximate the state $u_h$ of the \fom\ for an arbitrary parameter value $\theta' \in \Theta$. In general, $\theta'$ is an out-of-sample parameter, i.e., $\theta' \not \in \Theta_s$. The surrogate model constructed in \Cref{e:hosvd-stateapprox} only works for in-sample parameter values. We construct an interpolant in the parameter space to be able to evaluate the \htrom\ for out-of-sample $\theta'$. The associated \trom\ is referred to as \emph{interpolatory} TROM~\cite{mamonov2022:interpolatory}. We use conceptual ideas from~\cite{mamonov2022:interpolatory} to explain this approach. We first describe the in-sample case; the extension to out-of-sample $\theta'$ follows.

Suppose $\theta'$ belongs to $\Theta_s$. We define a vector $e^i(\theta) = (e_1^i(\theta), \ldots, e^i_{n_{\theta_i}}(\theta)) \in \{0,1\}^{n_{\theta_i}}$, $i = 1,\ldots,m$, as
\begin{equation}\label{e:extraction-vector}
e^i_j(\theta')=
\begin{cases}
1 & \text{if}\,\,\theta_i' = \theta_i^j\\
0 & \text{otherwise}
\end{cases}
\end{equation}

\noindent for $j = 1, \ldots, n_{\theta_i}$, where $\theta_i^j$ is the $j$th sample for the $i$th entry of the $m$-dimensional parameter vector within the interval $\Theta_i = [\theta_i^{\text{min}}, \theta_i^{\text{max}}]$. The vector $e^i$ encodes the position of $\theta_i'$ among the grid nodes in $\Theta_i$. Using this construction, we can extract snapshots corresponding to a particular $\theta' \in \Theta_s$ based on the operation
\begin{equation}\label{e:snapshot-extraction}
S(\theta') = \mathcal{S} \times_3 e^1(\theta') \times_4 e^2(\theta') \cdots \times_{m + 2} e^m(\theta') \in \mathbb{R}^{n_x \times n_t},
\end{equation}

\noindent where $\times_k$ denotes the $k$-mode tensor-vector product; let $\mathcal{X} \in \mathbb{R}^{m_1 \times \ldots \times m_l}$ be a tensor of order $l$ and $y \in \mathbb{R}^{m_k}$, then $\mathcal{X} \times_k y$ results in a tensor of order $l-1$ and size $m_1 \times \cdots \times m_{k-1} \times m_{k+1} \times \cdots \times m_l$. The operation in \Cref{e:snapshot-extraction} extracts snapshots for a particular $\theta' \in \Theta_s$ from the snapshot tensor $\mathcal{S}$ (i.e., it implements \Cref{e:hosvd-stateapprox} for all $j = 1, \ldots, n_t$ and $x_k \in \Omega^h$).

To construct the \trom, we introduce an interpolation procedure $\ell^i : \Theta \to \mathbb{R}^{n_{\theta_i}}$, $i = 1,\ldots,m$. The $j$th entry $\ell^i_j$ is defined by a Lagrangian polynomial of order $p$. We are given $\theta' \in \Theta$. Let $\theta_i^{i_1}, \ldots, \theta_i^{i_{p+1}}$ denote the $p+1$ closest grid nodes to $\theta_i'$ on $\Theta_i$. The $j$th entry of $\ell^i$ is given by
\begin{equation}\label{e:interpolation-polynomial}
\ell^i_j(\theta') =
\begin{cases}
\prod_{l=1, l \not=k}^{p+1} (\theta_i^{i_l} - \theta_i') / \prod_{l=1, l \not=k}^{p+1} (\theta_i^{i_l} - \theta_i^j) & \text{if}\,\, j = i_k \in \{i_1,\ldots,i_{p+1}\}, \\
0 & \text{otherwise}.
\end{cases}
\end{equation}

For the numerical experiments reported in this study, we consider polynomial order $p=3$. The vectors $\ell^i \in \mathbb{R}^{n_{\theta_i}}$ extend the notion defined in \Cref{e:extraction-vector} to out-of-sample parameter vectors $\theta' \not \in \Theta_s$. This allows us to generalize \Cref{e:snapshot-extraction} to out-of-sample parameter vectors $\theta'$; we obtain
\begin{equation}\label{e:snapshot-interpolation}
\tilde{S}(\theta') = \tilde{\mathcal{S}} \times_3 \ell^1(\theta') \times_4 \ell^2(\theta') \cdots \times_{m + 2} \ell^m(\theta') \in \mathbb{R}^{n_x \times n_t}
\end{equation}

\noindent for the low-rank tensorial approximation $\tilde{\mathcal{S}}$ defined in \Cref{e:hosvd-lr-snapshot} of the snapshot tensor $\mathcal{S}$.

For a fixed $\theta'$, we obtain the parameter-specific reduced basis $\Psi$
from the first $r$ singular vectors of $\tilde{S}(\theta')$. In practice, we do not need to compute the truncated SVD of $\tilde{S} \in \mathbb{R}^{n_x \times n_t}$. We can construct the reduced basis $\Psi \in \mathbb{R}^{n_x \times r}$ using a truncated SVD of the small-size parameter-specific tensor core $\tilde{C}(\theta) \in \mathbb{R}^{r_x \times r_t}$~\cite{mamonov2022:interpolatory}. This is due to the fact that we can express $\tilde{S}(\theta')$ as
\begin{equation}\label{e:snapshot-reduced-space-approx}
\tilde{S}(\theta) = \Phi^{(x)} \tilde{C}(\theta) (\Phi^{(t)})^\mathsf{T} \approx (\Phi^{(x)}U_r)\Sigma_r (\Phi^{(t)}V_r)^\mathsf{T},
\end{equation}

\noindent with $\tilde{C}(\theta)\approx  U_r\Sigma_rV_r^\mathsf{T}$, $U_r \in \mathbb{R}^{r_x \times r}$, $\Sigma_r \in \mathbb{R}^{r \times r}$, $V_r \in \mathbb{R}^{r_t \times r}$ with user-defined target rank $r$. The parameter-specific tensor core is computed as
\begin{equation}\label{e:core-tensor-contraction}
\tilde{C}(\theta) = \mathcal{C} \times_3 (\Phi^{(\theta_1)}\ell^1(\theta)) \cdots \times_{m + 2} (\Phi^{(\theta_m)}\ell^m(\theta)) \in \mathbb{R}^{r_x \times r_t}.
\end{equation}

\noindent We refer to this step as \emph{tensor core contraction}. This contraction step yields a matrix that couples the spatial and temporal reduced dimensions. The truncated left-singular vectors of the parameter-specific tensor core $\tilde{C}(\theta)$ provide us with coordinates of the parameter-specific reduced basis in $\Phi^{(x)}$, i.e., $\Psi =  \Phi^{(x)} U_r \in \mathbb{R}^{n_x \times r}$.

\emph{Intrusive TROM}: Using this basis, we can---akin to the POD approach in \Cref{s:pod}---solve \Cref{e:semi-discretized-ode} by solving the reduced system
\begin{equation}\label{e:reduced-system-trom}
\text{d}_t \mu_h(t;\theta)
= f_h(\mu_h(t; \theta), \theta, t)
\end{equation}

\noindent for $\mu_h \in \mathbb{R}^r$. The corresponding physical state is given by $u_h(t;\theta) =\Psi \mu_h(t;\theta)$ for all $t$. We refer to this mode of operation of the interpolatory \htrom\ as \emph{intrusive}.

\emph{Non-intrusive TROM}: As an alternative, we also consider a \emph{non-intrusive} variant that avoids integrating \Cref{e:reduced-system-trom} in time (projection-based \htrom). Instead, we compute the truncated SVD of $\tilde{C}(\theta)$ to obtain the low rank approximation $\tilde{C}(\theta) \approx [U_r\Sigma_r V_r^\mathsf{T}](\theta)$ and subsequently project to the full space via the left and right action of $\Phi^{(x)}$ and $\Phi^{(t)}$, respectively; that is, $(\Phi^{(x)} U_r(\theta))\Sigma_r(\theta) (\Phi^{(t)} V_r(\theta))^\mathsf{T}$ as in \Cref{e:snapshot-reduced-space-approx}.

We summarize the \emph{offline stage} for constructing and evaluating the interpolatory \htrom\ in \Cref{a:hosvd-trom-offline}. We show the steps required during the \emph{online stage} for the intrusive variant of the interpolatory \htrom\ in \Cref{a:hosvd-trom-online}. The non-intrusive variant is summarized in \Cref{a:hosvd-trom-online-ni}. The offline stage is the most expensive part of the interpolatory \htrom; it requires repeated high-fidelity simulations by invoking the \fom\ for $n_{\theta} = \prod_{i=1}^m n_{\theta_i}$ parameter samples $\theta^l \in \Theta_s$. However, it needs to be performed only once.

\subsubsection{Interpolatory Tensor Train Decomposition}\label{s:tt}

The second low-rank tensor decomposition we consider is the Tensor Train (\textbf{TT}) decomposition. TT decompositions provide a low-rank representation of high-order tensors with storage that scales linearly in the tensor order~\cite{oseledets2011:tensor_train}. In contrast to the Tucker representation in \Cref{s:tucker}, which compresses each mode by a single global factor matrix, TT represents a tensor through a chain of three-way cores with moderate intermediate ranks.

Let $\mathcal{S} \in \mathbb{R}^{n_x\times n_t\times n_{\theta_1}\times \cdots\times n_{\theta_m}}$ again denote the order-$(m + 2)$ tensor representation of the snapshot matrix $S$. The TT approximation $\tilde{\mathcal{S}}\in\mathbb{R}^{n_x\times n_t\times n_{\theta_1}\times\cdots\times n_{\theta_m}}$ of $\mathcal{S}$ is defined entry-wise: its $(i_1,\ldots,i_{m+2})$th scalar entry is given by
\[
\tilde{\mathcal{S}}(i_1,\ldots,i_{m+2})
=
\mathcal{G}^{(1)}(\,:\,,i_1,\,:\,)\,\mathcal{G}^{(2)}(\,:\,,i_2,\,:\,)\, \cdots\, \mathcal{G}^{(m+2)}(:,i_{m+2},:),
\]

\noindent where $\mathcal{G}^{(k)}\in\mathbb{R}^{r_{k-1}\times n_k\times r_k}$, $k = 1,\ldots,m+2$, are the TT cores; $n_1=n_x$, $n_2=n_t$, $n_{k+2}=n_{\theta_k}$ are introduced for notational convenience; and $r_0=r_{m+2}=1$ by construction.

\paragraph{Offline Stage} We start by describing the \emph{offline stage} of the \ttrom. As in \Cref{s:tucker}, we represent the snapshot matrix $S$ as a tensor $\mathcal{S}$ of size $n_x \times n_t \times n_{\theta_1} \times \cdots \times n_{\theta_m}$ and order $m+2$. We compute a TT approximation $\tilde{\mathcal{S}}$ of $\mathcal{S}$ using the TT-SVD procedure~\cite{oseledets2011:tensor_train}, i.e., a sequence of truncated SVDs, with truncation ranks selected according to \Cref{e:energy-bound-r} with tolerance $\epsilon^{\mathit{trom}}_{\mathit{off}} \in (0,1)$.

For notational convenience, we again set $n_1=n_x$, $n_2=n_t$, $n_{k+2}=n_{\theta_k}$ for $k=1,\ldots,m$, so that $\mathcal{S}$ is of size $n_1\times\cdots\times n_{m+2}$. The TT-SVD constructs TT cores $\mathcal{G}^{(k)}\in\mathbb{R}^{r_{k-1}\times n_k\times r_k}$ with $r_0=r_{m+2}=1$ by proceeding left to right. We initialize $\mathcal{Q}^{(1)}\coloneqq\mathcal{S}$. At stage $k=1,\ldots,m+1$, we view the current tensor $\mathcal{Q}^{(k)}$ as an element of $\mathbb{R}^{r_{k-1}\times n_k\times n_{k+1}\times\cdots\times n_{m+2}}$ and form the compound unfolding $Q^{(k)}$ that groups the first two modes into rows and the remaining modes into columns,
\begin{equation}\label{e:tt-unfold-offline}
Q^{(k)} \in \mathbb{R}^{(r_{k-1}n_k)\times(n_{k+1} \cdots n_{m+2})}.
\end{equation}

We compute a truncated SVD of $Q^{(k)}$ in \Cref{e:tt-unfold-offline}, selecting $r_k$ based on the criterion \Cref{e:energy-bound-r} with tolerance $\epsilon^{\mathit{trom}}_{\mathit{off}} \in (0,1)$, to obtain
\begin{equation}\label{e:tt-svd-stage-k}
Q^{(k)} \;\approx\;
U^{(k)}_{r_k}\,\Sigma^{(k)}_{r_k}\,\bigl(V^{(k)}_{r_k}\bigr)^{\mathsf{T}},
\qquad U^{(k)}_{r_k}\in\mathbb{R}^{(r_{k-1}n_k)\times r_k}.
\end{equation}

The $k$th TT core $\mathcal{G}^{(k)}$ is obtained by mapping $U^{(k)}_{r_k}$ into a three-way tensor format of size $r_{k-1} \times n_k \times r_k$; that is,
\begin{equation}\label{e:tt-core-from-U}
\mathcal{G}^{(k)} \;\gets\; \texttt{reshape}(U^{(k)}_{r_k}, (r_{k-1}, n_k, r_k)),
\qquad k=1,\ldots,m+1.
\end{equation}

To proceed to the next stage, we absorb the singular values and right singular vectors into the updated tensor $\mathcal{Q}^{(k+1)}$:
\begin{equation}\label{e:tt-update-next}
\mathcal{Q}^{(k+1)} \;\gets\; \texttt{reshape}\bigl(\Sigma^{(k)}_{r_k}\,\bigl(V^{(k)}_{r_k}\bigr)^{\mathsf{T}},\,(r_k, n_{k+1}, \ldots, n_{m+2})\bigr), \qquad k=1,\ldots,m+1.
\end{equation}

After completing stages $k=1,\ldots,m+1$, the remaining tensor $\mathcal{Q}^{(m+2)}$ has size $r_{m+1}\times n_{m+2}$. We define the final core by adding a singleton third mode:
\begin{equation}\label{e:tt-final-core}
\mathcal{G}^{(m+2)} \;\gets\; \texttt{reshape}\bigl(\mathcal{Q}^{(m+2)},\,(r_{m+1}, n_{m+2}, 1)\bigr).
\end{equation}

This enforces $r_{m+2}=1$ and completes the TT representation $\tilde{\mathcal{S}}$ through the chain of cores $\{\mathcal{G}^{(k)}\}_{k=1}^{m+2}$.

\paragraph{Online Stage} The \emph{online stage} of the \ttrom\ aims to rapidly approximate the state $u_h$ of the \fom\ for an arbitrary parameter value $\theta' \in \Theta$. In general, $\theta'$ is an out-of-sample parameter, i.e., $\theta' \not\in \Theta_s$. To enable evaluation at such $\theta'$, we use the same interpolation vectors $\ell^i(\theta') \in \mathbb{R}^{n_{\theta_i}}$ from \eqref{e:interpolation-polynomial} as in \Cref{s:tucker}, $i=1,\ldots,m$ (with $\ell^i=e^i$ for in-sample $\theta'\in\Theta_s$).

Recall the TT representation of $\tilde{\mathcal{S}}$ via cores $\mathcal{G}^{(k)}\in\mathbb{R}^{r_{k-1}\times n_k\times r_k}$, $k=1,\ldots,m+2$, where
\[
n_1=n_x,\qquad n_2=n_t,\qquad n_{l+2}=n_{\theta_l}\quad (l=1,\ldots,m),
\qquad r_0=r_{m+2}=1.
\]

\noindent Since $r_0=1$, the first core $\mathcal{G}^{(1)}\in\mathbb{R}^{1\times n_x\times r_1}$ has a singleton first mode; its mode-$2$ unfolding defines the spatial factor, obtained by collapsing the singleton first dimension:
\begin{equation}\label{e:tt-phi-x}
\Phi^{(x)} \gets \texttt{reshape}\bigl(\mathcal{G}^{(1)},\,(n_x, r_1)\bigr) \in \mathbb{R}^{n_x\times r_1}.
\end{equation}

To construct a parameter-specific matrix that couples the spatial rank and time, we contract the parameter cores using the interpolation vectors. For $i=1,\ldots,m$, let $k=i+2$ and define
\begin{equation}\label{e:tt-param-core-contraction}
M^{(k)}(\theta') = \sum_{j=1}^{n_{\theta_i}} \ell^i_j(\theta')\, \mathcal{G}^{(k)}(\,:\,,j,\,:\,) \in \mathbb{R}^{r_{k-1}\times r_k}.
\end{equation}

The product of these small matrices yields a vector
\begin{equation}\label{e:tt-v-theta}
v(\theta') = M^{(3)}(\theta')\,M^{(4)}(\theta')\cdots M^{(m+2)}(\theta')
\in \mathbb{R}^{r_2}.
\end{equation}

Using $v(\theta')$, we contract the time core $\mathcal{G}^{(2)}\in\mathbb{R}^{r_1\times n_t\times r_2}$ along its third mode and obtain the parameter-specific matrix
\begin{equation}\label{e:tt-core-tensor-contraction}
\tilde{C}(\theta') = \mathcal{G}^{(2)} \times_3 v(\theta') \in \mathbb{R}^{r_1\times n_t}.
\end{equation}

Finally, the TT approximation of the snapshot matrix at $\theta'$ can be written as
\begin{equation}\label{e:tt-snapshot-reduced-space-approx}
\tilde{S}(\theta') \;=\; \Phi^{(x)}\,\tilde{C}(\theta') \in \mathbb{R}^{n_x\times n_t}.
\end{equation}

For a fixed $\theta'$, we obtain a parameter-specific reduced basis from the truncated SVD of $\tilde{C}(\theta')$,
\begin{equation}\label{e:tt-svd-tildeC}
\tilde{C}(\theta') \approx U_r(\theta')\,\Sigma_r(\theta')\,\bigl(V_r(\theta')\bigr)^\mathsf{T},
\qquad U_r(\theta')\in\mathbb{R}^{r_1\times r},
\end{equation}

\noindent with $r(\epsilon^{\mathit{trom}}_{\mathit{on}})$ chosen via \Cref{e:energy-bound-r}. We define the full-space basis $\Phi^{(x)}U_r(\theta')\in\mathbb{R}^{n_x\times r}$. Using this basis, we can---akin to the POD approach in \Cref{s:pod}---solve \Cref{e:semi-discretized-ode} by integrating the reduced system \Cref{e:reduced-system-trom}. We refer to this mode of operation of the interpolatory \ttrom\ as \emph{intrusive}. As an alternative, a \emph{non-intrusive} variant avoids time integration by reconstructing
\[
\hat{X}(\theta') = \bigl(\Phi^{(x)}U_r(\theta')\bigr)\Sigma_r(\theta')V_r(\theta')^\mathsf{T}\in\mathbb{R}^{n_x\times n_t},
\]

\noindent and taking its columns as time snapshots.

We summarize the offline and online stages of the \ttrom\ in \Cref{a:tt-trom-offline,a:tt-trom-online,a:tt-trom-online-ni}, respectively.

\subsubsection{Pre-Compression in Space}\label{s:precompression}

For problems with large spatial dimension $n_x$ and high-dimensional parameter spaces, constructing and storing the full snapshot tensor $\mathcal{S}$ of size $n_x \times n_t \times n_{\theta_1}\times\cdots\times n_{\theta_m}$ may become prohibitive. In particular, all three offline constructions in \Cref{s:pod}--\Cref{s:tt} involve (explicitly or implicitly) the mode-$1$ unfolding resulting in a matrix of size $n_x \times n_t n_{\theta}$ with $ n_{\theta} \coloneqq \prod_{i=1}^m n_{\theta_i}$, whose leading dimension is the full spatial size $n_x = \prod_{i=1}^d n_i$. If we were to carry out the offline construction in the full space, then explicitly forming and manipulating this matrix would require memory and offline wall-clock times that are intractable even on many modern high-performance computing platforms. To alleviate this bottleneck, we apply a preliminary \emph{pre-compression} in the spatial direction. The main idea is to first identify a low-dimensional spatial subspace that captures the dominant spatial features of the solution manifold, and then project all snapshots onto this subspace before constructing the \rom\ surrogates.

\paragraph{Offline stage}
Let $\Theta_s=\{\theta^l\}_{l\in\mathcal{I}_s}$ denote the full training set used for snapshot generation, as in \Cref{s:pod,s:tucker,s:tt}. We select a (typically much smaller) subset $\Theta_{\mathit{pre}}=\{\theta^{l}\}_{l\in\mathcal{I}_{\mathit{pre}}}\subset \Theta_s$, $n_{\theta,\mathit{pre}} = |\mathcal{I}_{\mathit{pre}}| \ll n_{\theta}$. We form the snapshot matrix
\begin{equation}\label{e:precompression-matrix-new}
S_{\mathit{pre}}
=
\bigl[
u_h^1(\theta^{1,\ldots,1}),u_h^2(\theta^{1,\ldots,1}),
\ldots,
u_h^{n_t}(\theta^{n_{\theta_1,\mathit{pre}},\ldots,n_{\theta_m,\mathit{pre}}})
\bigr]
\in \mathbb{R}^{n_x \times n_t n_{\theta,\mathit{pre}}}.
\end{equation}

\noindent We compute a truncated (thin) SVD of $S_{\mathit{pre}}$. We choose the truncation rank $r_x$ based on \Cref{e:energy-bound-r} with tolerance $\epsilon^{\mathit{pre}}\in(0,1)$ to obtain $S_{\mathit{pre}} \approx U_{\mathit{pre}}\Sigma_{\mathit{pre}} V_{\mathit{pre}}^\mathsf{T}$. We set $\Phi^{(x)}_{\mathit{pre}} = U_{\mathit{pre}}\in\mathbb{R}^{n_x\times r_x}$ with $(\Phi^{(x)}_{\mathrm{pre}})^\mathsf{T}\Phi^{(x)}_{\mathrm{pre}}=I_{r_x}$; the columns of $\Phi^{(x)}_{\mathit{pre}}$ define the pre-compressed spatial basis.

For each training parameter $\theta^l\in\Theta_s$ and time index $j=1,\ldots,n_t$, we project the corresponding \fom\ snapshot onto the pre-compressed space via
\begin{equation}\label{e:precompression-projection-new}
\hat{u}^j_h(\theta^l)
=
(\Phi^{(x)}_{\mathit{pre}})^\mathsf{T} u_h^j(\theta^l)
\in \mathbb{R}^{r_x}.
\end{equation}

\noindent Assembling these reduced snapshots yields the compressed $r_x \times n_t \times n_{\theta_1}\times\cdots\times n_{\theta_m}$ snapshot tensor $\hat{\mathcal{S}}$ according to
\[
\hat{\mathcal{S}}(\,:\,,j,l_1,\ldots,l_m) = \hat{u}_h^j(\theta^l).
\]

\noindent The tensor $\hat{\mathcal{S}}$ replaces $\mathcal{S}$ in all subsequent offline stages.

\begin{itemize}[leftmargin=*]
\item
\textbf{\podrom:} Form the snapshot matrix $\hat{S}\in\mathbb{R}^{r_x\times n_t n_\theta}$ from the tensor $\hat{\mathcal{S}}$ and compute the associated POD basis $\hat{\Phi}\in\mathbb{R}^{r_x\times r}$ using the tolerance $\epsilon^{\mathit{pod}}$. The lifted full-space basis is given by
\[
\Phi = \Phi^{(x)}_{\mathrm{pre}}\,\hat{\Phi}\in\mathbb{R}^{n_x\times r}.
\]
\item
\textbf{\htrom:} We apply the HOSVD to $\hat{\mathcal{S}}$ to obtain $\hat{\Phi}^{(x)}\in\mathbb{R}^{r_x\times r_x'}$, $\Phi^{(t)}$, $\{\Phi^{(\theta_i)}\}_{i=1}^m$ and a core tensor $\hat{\mathcal{C}}$, with ranks selected via \Cref{e:energy-bound-r} using $\epsilon^{\mathit{trom}}_{\mathit{off}}$. The lifted \emph{full-space} basis is given by
\[
\Phi^{(x)} = \Phi^{(x)}_{\mathit{pre}}\hat{\Phi}^{(x)} \in \mathbb{R}^{n_x\times r_x'}.
\]
\item
\textbf{\ttrom:} We apply the TT-SVD to $\hat{\mathcal{S}}$ (with leading mode size $r_x$) to obtain TT cores $\{\hat{\mathcal{G}}^{(k)}\}_{k=1}^{m+2}$ with ranks chosen according to \Cref{e:energy-bound-r} using $\epsilon^{\mathit{trom}}_{\mathit{off}}$. In the TT online stage, the spatial factor extracted from the first core $\hat{\Phi}^{(x)}=\texttt{reshape}(\hat{\mathcal{G}}^{(1)},\,(r_x, r_1))\in\mathbb{R}^{r_x\times r_1}$ is lifted to the full space via
\[
\Phi^{(x)} = \Phi^{(x)}_{\mathit{pre}}\,\hat{\Phi}^{(x)} \in \mathbb{R}^{n_x\times r_1}.
\]
\end{itemize}

\paragraph{Online stage.} All online evaluations (intrusive or non-intrusive) are performed in the pre-compressed space of dimension $r_x$. Let $\hat{u}_h^j(\theta)\in\mathbb{R}^{r_x}$ denote the resulting reduced-space approximation. The corresponding approximation is lifted to the original space via the map
\[
u_h^j(\theta) \approx \Phi^{(x)}_{\mathit{pre}}\,\hat{u}_h^j(\theta)\in\mathbb{R}^{n_x}.
\]

\noindent Equivalently, whenever a method produces a parameter-specific basis acting on the pre-compressed coordinates (e.g., $\hat{\Phi}\Psi(\theta)$ in \Cref{s:tucker,s:tt}), the associated full-space basis is obtained by left multiplication with $\Phi^{(x)}_{\mathit{pre}}$.

We summarize the pre-compression in \Cref{a:precompression-offline}.

\section{Numerical Results}\label{s:results}

Next, we present empirical results for the performance of the considered {\rom}s and explore the behavior of the developed methodology to changes in the hyperparameter choices.

\subsection{Performance Measures} At any snapshot time $t^j \in [0,1]$, $j = 1,\dots,n_t$, the approximate full-order state is reconstructed by projecting the solution of the reduced system to the full space. We denote the high-fidelity \fom\ solution at final time $t = 1$ by $u_{\mathit{fom}}$ and the corresponding \rom\ solution by $u_{\mathit{rom}}$. We compute the reconstruction error $\varepsilon_p$ based on the relative $\ell^2$- and $\ell^\infty$-norms of the residual $u_{\mathit{fom}} - u_{\mathit{rom}}$. That is,
\[
\varepsilon_p = \|u_{\mathit{fom}} - u_{\mathit{rom}}\|_p / \|u_{\mathit{fom}}(t^j)\|_p, \qquad p = 2, \infty.
\]

\noindent In addition to reconstruction errors, we also report run times for the different \rom\ variants and the \fom.

\subsection{Algorithmic Parameters}

We select the optimal ranks for all \rom\ variants according to \Cref{e:energy-bound-r}. In addition, we provide results in which we match the \podrom\ offline rank with the \trom\ rank determined during its online phase. This allows us to compare the \podrom\ to the \trom\ variants using two regimes, one in which we construct an accurate representation (computing the optimal \podrom\ rank according to \Cref{e:energy-bound-r}) and one in which we match the optimal compression level determined for the \trom.

We need to determine several algorithmic parameters that determine the performance of the designed numerical approach. We did so empirically. For the snapshot generation and solution of the problem, we need to determine the number of time steps $n_t \in \mathbb{N}$ (i.e., time step size $h_t$) that balances numerical accuracy against computational performance. We note that the considered numerical scheme is unconditionally stable, i.e., $n_t$ can be determined solely on accuracy requirements and not based on stability concerns (see \Cref{s:time_integration}). If not noted otherwise, we select $n_t = 16$.

Moreover, we need to determine the oversampling rate $p_o \in \mathbb{N}$ (Line~\ref{l:oversampling} in \Cref{a:rsvd-fixed-accuracy}) and the target rank increment $p_s\in \mathbb{N}$ (Line~\ref{l:range-finder} in \Cref{a:rsvd-fixed-accuracy}) for the rSVD. These rSVD parameters are fixed throughout all experiments. We choose $p_o = 16$ and $p_s = 32$.

Moreover, we need to determine the tolerances for the optimal rank selection using \Cref{e:energy-bound-r} during the online and offline stages of the different \rom\ variants considered here. We determined these tolerances empirically. To do so, we reduced them by one order of magnitude until a further reduction did not yield any significant changes in the reconstruction accuracy (the errors remained at the same order). Our choices are $\epsilon^{\mathit{trom}}_{\mathit{off}} = \text{1e--5}$ for the time mode,  $\epsilon^{\mathit{trom}}_{\mathit{off}} = \text{1e--6}$ for the space mode, and $\epsilon^{\mathit{trom}}_{\mathit{off}} = \text{1e--7}$ for the parameter modes for the offline stage for both \trom-unfoldings (\htrom\ and \ttrom). The tolerance for the \podrom\ is set to $\epsilon^{\mathit{pod}} = \text{1e--6}$. The tolerance for the {\trom}s for the online stage is $\epsilon^{\mathit{trom}}_{\mathit{on}} = \text{1e--6}$. The tolerance for the pre-compression in space is $\epsilon^{\mathit{pre}} =  \text{1e--8}$.

Lastly, we need to decide for which parameter values $\theta_i$ we generate snapshots. We draw $n_{\theta_i}$, $i = 1,\ldots,m$, equispaced samples from the hypercube $\Theta = \bigotimes_{i=1}^m [\theta_i^{\mathit{min}}, \theta_i^{\mathit{max}}] \subset \mathbb{R}_+^m$. We empirically select $n_{\theta_1} = 3$ (for $\alpha$), $n_{\theta_2} = 10$ (for $\rho$), $n_{\theta_3} = 3$ (for $u_o^{\mathit{inv}}$), $n_{\theta_4} = 3$ (for $u_o^{\mathit{hyp}}$), $n_{\theta_5} = 3$ (for $\gamma_p$), $n_{\theta_6} = 3$ (for $\gamma_i$), $n_{\theta_7} = 3$ (for $\lambda_d$), $n_{\theta_8} = 5$ (for $\kappa_s$), and $n_{\theta_9} = 5$ (for $\kappa_c$).  We report experiments below to justify these choices. The total number of snapshots is $n_{\theta} = \prod_{i=1}^9 n_{\theta_i} = \inum{182250}$.

We summarize the choices for the values of all algorithmic parameters in \Cref{t:rom-hyperparameters}.

During some of our experiments, we randomly draw parameter samples $\theta \in \mathbb{R}^m$. We only accept samples that satisfy the model constraint $\uo^{\textit{hyp}} \leq \uo^{\textit{inv}}$.

\begin{table}
\caption{Algorithmic parameters that control the performance of the considered {\rom}s. We report (from left to right), the symbol, the meaning, a reference to identify the context in which the parameter appears, and the default value considered in this study (determined empirically). We note that we report numerical experiments in \Cref{s:results} that justify the selection for some of these default values.}
\label{t:rom-hyperparameters}
\tabadjust
\begin{tabular}{llll}
\toprule
\bf symbol & \bf meaning & \bf reference & \bf recommended value\\
\midrule
$n_t\in\mathbb{N}$ & number of time steps for time integrator & \Cref{s:time_integration} & 16 \\
$p_o \in \mathbb{N}$ & oversampling (rSVD) & Line~\ref{l:oversampling} in \Cref{a:rsvd-fixed-accuracy} & 16  \\
$p_s \in \mathbb{N}$ & target rank increment (adaptive rSVD) & Line~\ref{l:range-finder} in \Cref{a:rsvd-fixed-accuracy} & 32\\
$p \in \mathbb{N}$ & polynomial interpolation order & \Cref{e:interpolation-polynomial} & 3  \\
$\epsilon^{\mathit{trom}}_{\mathit{off}} > 0$ & tolerance for \trom\ offline rank selection & \Cref{e:energy-bound-r} &
1e--5 (time) \\
& & & 1e--6 (space) \\
& & & 1e--7 (parameters) \\
$\epsilon^{\mathit{trom}}_{\mathit{on}} > 0$ & tolerance for \trom\ online rank selection  & \Cref{e:energy-bound-r} & 1e--6 \\
$\epsilon^{\mathit{pod}} > 0$ & tolerance for POD rank selection & \Cref{e:energy-bound-r} & 1e--6 \\
$\epsilon^{\mathit{pre}} > 0$ & tolerance for pre-compression of space & \Cref{s:precompression} & 1e--8 \\
$n_{\theta_i} \in \mathbb{N}$ & number of samples per parameter $\theta_i$ & \Cref{e:multiindexset-para} & see main text \\
\bottomrule
\end{tabular}
\end{table}

\subsection{Considered \rom\ Variants}

We select several \rom\ variants for comparison. The \podrom\ serves as a baseline. We consider two \trom\ variants: \htrom\ and \ttrom. For each \trom\ variant we consider intrusive and non-intrusive (projection based) approaches. Overall, we observed that the \htrom\ and \ttrom\ perform similarly. Consequently, we decided to perform most of the experiments focusing on the \htrom\ variant.

\subsection{Baseline Experiment: Evaluation of the \fom}

\ipoint{Purpose} This experiment serves as a baseline in which we explore performance of our solver for the \fom.

\ipoint{Setup} We execute the \fom\ solver for different mesh sizes (the native resolution of the three-dimensional data is $193 \times 229 \times 193$). We set the number of time steps of our solver to $n_t = 32$. We set the model parameters to the default values listed in \Cref{t:multispecies-para}.

\begin{figure}
\includegraphics[width=0.5\textwidth]{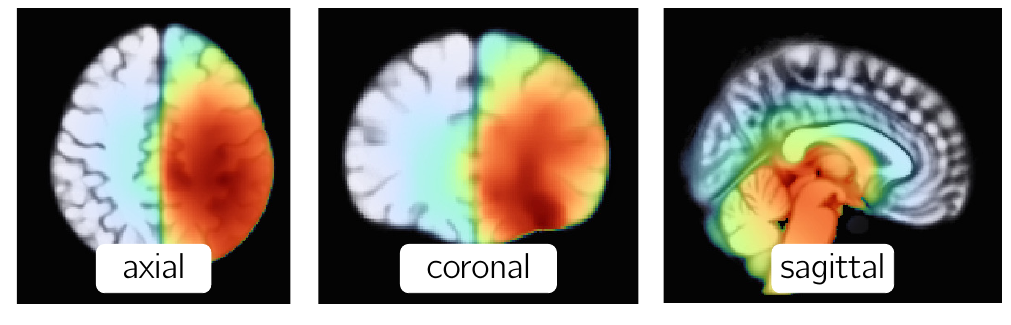}
\caption{Exemplary simulation results for a single-species model in three dimensions (ambient space). We solve the \fom. We show (from left to right) an axial, coronal and sagittal view of the 3D volume. The native resolution (number of mesh points) of this data set is $193 \times 229 \times 193$ at pseudo-time $t=1$. High tumor cell densities are shown in red and low densities are shown in green.}
\label{f:fom-3D}
\end{figure}

\begin{figure}
\includegraphics[width=0.8\textwidth]{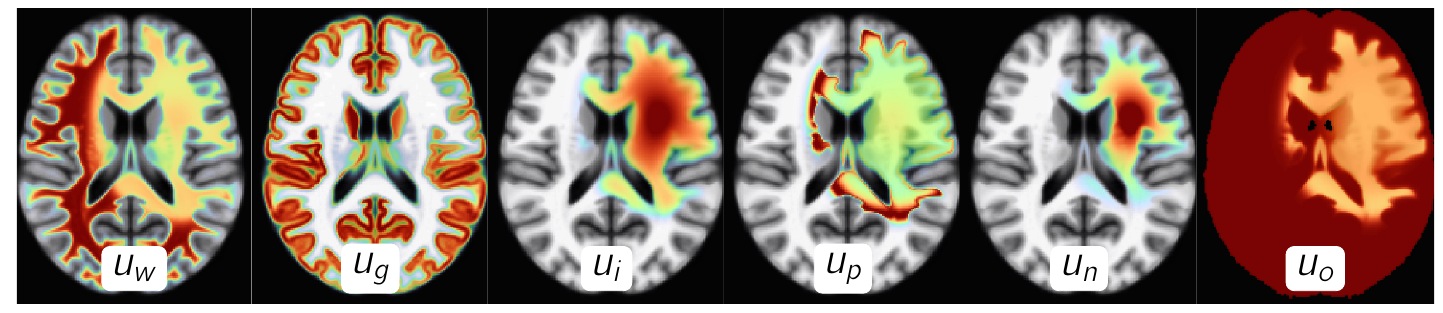}
\caption{Simulation results for the multi-species model in two dimensions (ambient space). We solve the \fom. We show (from left to right) axial slices for the white matter density $\uw$, the gray matter density $\ug$, and the densities for the infiltrative tumor cells $\ui$, proliferative tumor cells $\up$, and necrotic tumor cells $\un$, as well as the oxygen concentration $\uo$ at pseudo-time $t=1$. High cell densities are shown in red and low densities are shown in green.}
\label{f:fom-2D}
\end{figure}

\ipoint{Results} We show exemplary results for a three-dimensional simulation (ambient space) for the single-species model in \Cref{f:fom-3D}. We show results for a two-dimensional simulation (ambient space) for the multi-species model in \Cref{f:fom-2D}. We report runtime results in \Cref{t:runtime-fom-3D} and \Cref{t:runtime-fom-2D-ms}, respectively, as a function of the mesh size.

\begin{table}
\caption{Runtime for our solver for the \fom\ for the three-dimensional single-species tumor growth model. We consider three mesh sizes (the bottom row is for the original resolution). We report (from left to right) the mesh size, and runtimes (in seconds) for the two diffusion steps and the reaction step in our splitting approach. The last column shows the total runtime.}
\label{t:runtime-fom-3D}
\tabadjust
\begin{tabular}{crrrr}\toprule
\bf mesh  & \multicolumn{2}{c}{\bf diffusion} & \bf reaction & \bf runtime  \\
& \bf step 1 & \bf step 2 \\
\midrule
$100\times 125 \times 100$ & \fnum{ 14.795} & \fnum{ 17.078} & \fnum{0.502089} & \fnum{  32.71000} \\
$200\times 250 \times 200$ & \fnum{140.826} & \fnum{156.849} & \fnum{3.842000} & \fnum{ 304.24800} \\
\midrule
$193\times 229 \times 193$ & \fnum{127.495} & \fnum{140.726} & \fnum{3.318000} & \fnum{273.919875} \\
\bottomrule
\end{tabular}
\end{table}

\begin{table}
\caption{Runtime(s) for the FOM for the 2D multispecies tumor growth model for three different mesh sizes (the bottom row is for the original resolution). We report (from left to right) the mesh size and runtimes (in seconds) for the diffusion steps and the reaction step in our splitting approach. The last column shows the total runtime.}
\label{t:runtime-fom-2D-ms}
\tabadjust
\begin{tabular}{crrrr}\toprule
\bf mesh  & \multicolumn{2}{c}{\bf diffusion} & \bf reaction & \bf runtime  \\
& \bf step 1 & \bf step 2 \\
\midrule
$100\times 125$ & \fnum{0.175} & \fnum{0.188} & \fnum{0.323} & \fnum{1.576} \\
$200\times 250$ & \fnum{0.970} & \fnum{0.937} & \fnum{0.656} & \fnum{5.533} \\
\midrule
$193\times 229$ & \fnum{0.506} & \fnum{0.478} & \fnum{0.410} & \fnum{2.170} \\
\bottomrule
\end{tabular}
\end{table}

\ipoint{Observations} Solving the problem using the \fom\ is expensive, particularly in three dimensions. This makes the \fom-solver impractical in many-query settings. We also observe that runtime grows linearly with mesh refinement, suggesting that our \fom-solver has good computational scaling (in space).

\subsection{Model Parameter Sensitivities}\label{s:parameter-sensitivity}

\ipoint{Purpose} We study the sensitivity of the model output (states) with respect to changes in the parameters $\theta$. Our goal is to determine if there are any parameters $\theta_i$ that do not really affect the model state $u$.

\ipoint{Setup} We consider the multi-species model in two dimensions. The model is controlled by $m=9$ parameters. We execute the \fom\ solver on the original data resolution ($193 \times 229$). In \Cref{s:surrogate} we introduced the parameter ranges $\Theta_i = [\theta_i^{\mathit{min}}, \theta_i^{\mathit{max}}] \subset \mathbb{R}_+$ for each $\theta_i$. We compute a reference solution at the center of each of these intervals by choosing
\[
\theta^{\mathit{ref}} = (\theta^\text{ref}_1,\ldots,\theta^\text{ref}_m) \in \Theta_i,
\quad
\theta_i^{\mathit{ref}} = (\theta_i^{\mathit{min}} + \theta_i^{\mathit{max}}) / 2.
\]

\noindent We study how perturbations $\delta\theta \in \mathbb{R}$, $\theta^{\mathit{trial}} = \theta^{\mathit{ref}} + \delta\theta e_i$, $\theta^{\mathit{trial}} \in \Theta$, with unit vector $e_i \in \{0,1\}^m$, $(e_i)_j = 1$ if $i=j$ and $(e_i)_j = 0$ otherwise, affect the state of the system. We report the $\ell^2$-norm of the residual $\varepsilon_{ij} \in \mathbb{R}$ of the state variables computed for $\theta^{\mathit{ref}}$ and $\theta^{\mathit{trial}}$ at final time $t=1$; that is,
\[
\varepsilon_{ij} = \|[u_h^{n_t}]_j(\theta^{\text{ref}}) - [u_h^{n_t}]_j(\theta^{\text{ref}} + \delta\theta e_i)\|_2^2,
\]

\noindent for the $j$th species $[u_h^{n_t}]_j \in \mathbb{R}^{n_x}$,  $j = 1,\ldots,n_s$, and the $i$th parameter $\theta_i$, $i=1,\ldots,m$. The perturbations $\delta\theta$ are selected so that we draw equispaced samples within each interval $\Theta_i$, spanning the entire interval.

\begin{figure}
\centering
\includegraphics[width=0.85\textwidth]{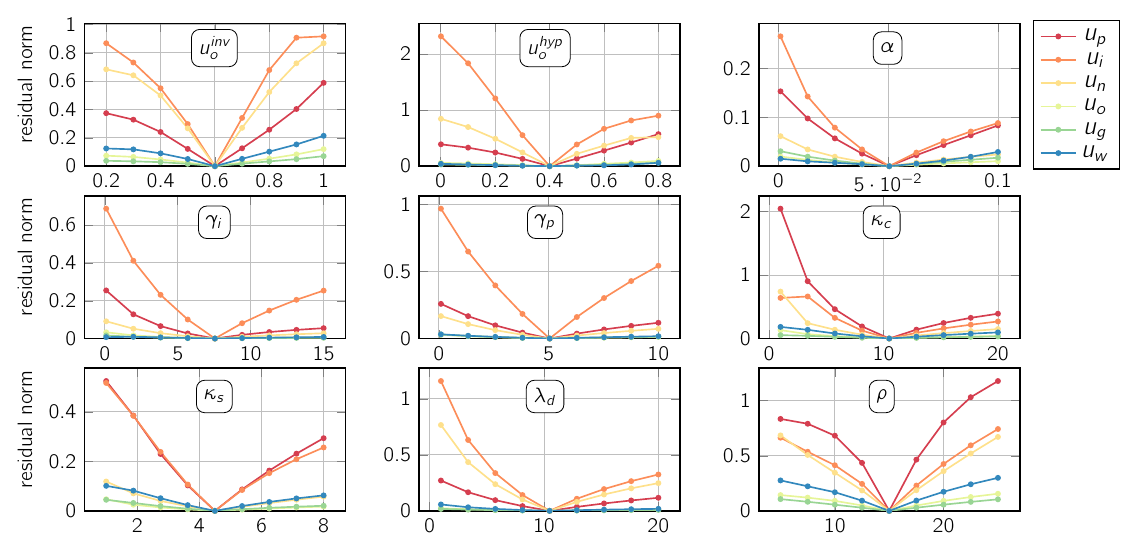}
\caption{Parameter sensitivities for the multi-species model. We plot the $\ell^2$-norm of the residual for each species ($\up$, $\ui$, $\un$, $\uo$, $\ug$, $\uw$) with respect to a reference solution located at the center of each parameter interval $\Theta_i = [\theta_i^{\mathit{min}}, \theta_i^{\mathit{max}}] \subset \mathbb{R}$ as a function of the perturbation from this reference parameter. The reference solution is computed at the center of the hypercube $\Theta$ defined by the parameter ranges $\Theta_i$. We report results for all nine model parameters $\theta = (\rho, \alpha, \gamma_p, \gamma_i,\lambda_d, \kappa_c, \kappa_s, u_o^\textit{hyp}, u_o^\textit{inv})$. The error is zero at the center of each parameter domain. As we move away from the center we expect the error to increase if the model is sensitive with respect to the considered parameter. Each plot provides the norm of the residuals for these perturbations for each individual parameter for all species (states). Averages are reported in \Cref{t:sensitivity-analysis-summary}.}
\label{f:sensitivity-analysis-summary}
\end{figure}

\ipoint{Results} We report average values for the norm of the residuals in \Cref{t:sensitivity-analysis-summary}. We show the trend of the norm of the residuals as a function of the perturbation from $\theta^{\text{ref}}$ in \Cref{f:sensitivity-analysis-summary}.

\begin{table}
\caption{Summary of sensitivity analysis for the multi-species model. We perturb each parameter $\theta_i$ from a reference solution sitting at the center of each parameter interval $\Theta_i = [\theta_i^{\mathit{min}}, \theta_i^{\mathit{max}}] \subset \mathbb{R}$. We consider 9 parameter samples within each interval $\Theta_i$. We report results for all nine model parameters $\theta = (\rho, \alpha, \gamma_p, \gamma_i,\lambda_d, \kappa_c, \kappa_s, u_o^\textit{hyp}, u_o^\textit{inv})$. We report the average and standard deviation across all species (model states) with respect to each model parameter $\theta_i$.}
\label{t:sensitivity-analysis-summary}
\tabadjust
\begin{tabular}{lccccccccc}\toprule
\bf statistic & $\rho$ & $\alpha$ & $\gamma_p$ & $\gamma_i$ & $\lambda_d$ & $\kappa_c$ & $\kappa_s$ & $\uo^{\textit{hyp}}$ & $\uo^{\textit{inv}}$ \\
\midrule
\bf mean & \fnum{0.314100760} & \fnum{0.040245530} & \fnum{0.093980524} & \fnum{0.054233353} & \fnum{0.121473750} & \fnum{0.189617450} & \fnum{0.093663896} & \fnum{0.291682170} & \fnum{0.252844980} \\
\bf std  & \fnum{0.303229630} & \fnum{0.057165211} & \fnum{0.182567490} & \fnum{0.114177800} & \fnum{0.217961810} & \fnum{0.334027680} & \fnum{0.126880770} & \fnum{0.461259130} & \fnum{0.281891440} \\
\bottomrule
\end{tabular}
\end{table}

\ipoint{Observations} The most important observation is that at least one of the states is sensitive to changes in each of the model parameters. This indicates that we have to construct the \rom{s} for all parameter modes. The states $\up$ and $\ui$ are most sensitive. The states $\uw$ and $\ug$ (healthy tissue) are least sensitive. The trends also reveal that as we perturb the parameters, there is no clear indication that we have to increase the sampling in a particular zone of the intervals $\Theta_i$; the trend is almost symmetric with respect to positive or negative perturbations in the parameter space. Consequently, we stipulate that equispaced samples drawn inside the hypercube $\Theta$ are adequate for constructing the {\rom}s.

\subsection{\trom\ Parameter Sampling}

\ipoint{Purpose} We study how the number of samples $n_{\theta_i}$ drawn for each parameter $\theta_i$, $i=1,\ldots,9$, to build the snapshot tensor $\mathcal{S}$ during the offline stage of the \htrom\ affects the reconstruction accuracy. The goal of this experiment is to determine a policy for selecting the number of samples  $n_{\theta_i}$.

\begin{figure}
\centering
\includegraphics[width=0.75\textwidth]{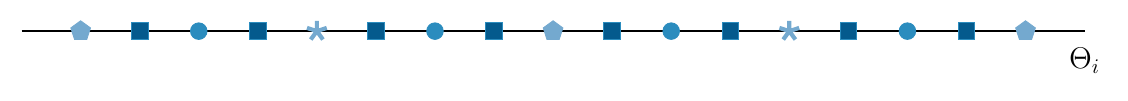}
\caption{Illustration of the refinement of the parameter samples in $\Theta_i$.  We gradually increase the number of parameter samples $n_{\theta_i}$ from 3 to  17. At each refinement step, we add new samples midway between the existing ones, while keeping all previously selected samples.}
\label{f:parameter-sample-refinement}
\end{figure}

\ipoint{Setup} We consider the multi-species model and place the parameter samples on an equispaced grid. {We build the \htrom\ for one parameter $\theta_i$ at a time. The remaining parameters $\theta_j$, $j \not=i$, remain fixed. The values for these $\theta_j$ are set to the center point of the interval $\Theta_j = [\theta_j^{\mathit{min}}, \theta_j^{\mathit{max}}$, i.e., $\theta_j = \theta_j^{\mathit{min}} + (\theta_j^{\mathit{max}} -\theta_j^{\mathit{min}})/2$ for all $j \not = i$. For the refinement of the sampling of $\Theta_i$ we} begin with $n_{\theta_i} = 3$ samples and gradually increase this number to  $n_{\theta_i} = 17$. At each refinement step, we add new samples midway between the existing ones, while keeping all previously selected samples. In this way, the parameter grid is successively refined. Accordingly, the \htrom\ is constructed using 3, 5, 9, and finally 17 parameter samples. This refinement process is illustrated in \Cref{f:parameter-sample-refinement}. We use $n_t = 16$ and the tolerances prescribed in \Cref{t:rom-hyperparameters}.

We additionally report results for $n_t = 64$ and $\epsilon^{\mathit{trom}}_{\mathit{on}}$, $\epsilon^{\mathit{trom}}_{\mathit{off}}$, and  $\epsilon^{\mathit{pre}}$ set to 1e--12. Changing only $n_t$ or the tolerances for the rank selection did not result in a significant reduction of the reconstruction error. We note that we also replaced the rSVD by a standard SVD for these high-accuracy experiments; we did not observe any deterioration of the performance; consequently, we use an rSVD. We also note that in practice we cannot afford using $n_t= 64$ when building the \htrom\ for all parameters due to runtime constraints and memory pressure.

After each refinement step we randomly draw 100 off-grid parameters $\theta_i^{\mathit{trial}} \in \Theta_i$, evaluate the \htrom, and compute the error between the \htrom\ solution and the \fom\ solution. We do this for each individual parameter. We solve the \fom\ model using $n_t = 16$.

\begin{figure}
\centering
\includegraphics[width=\textwidth]{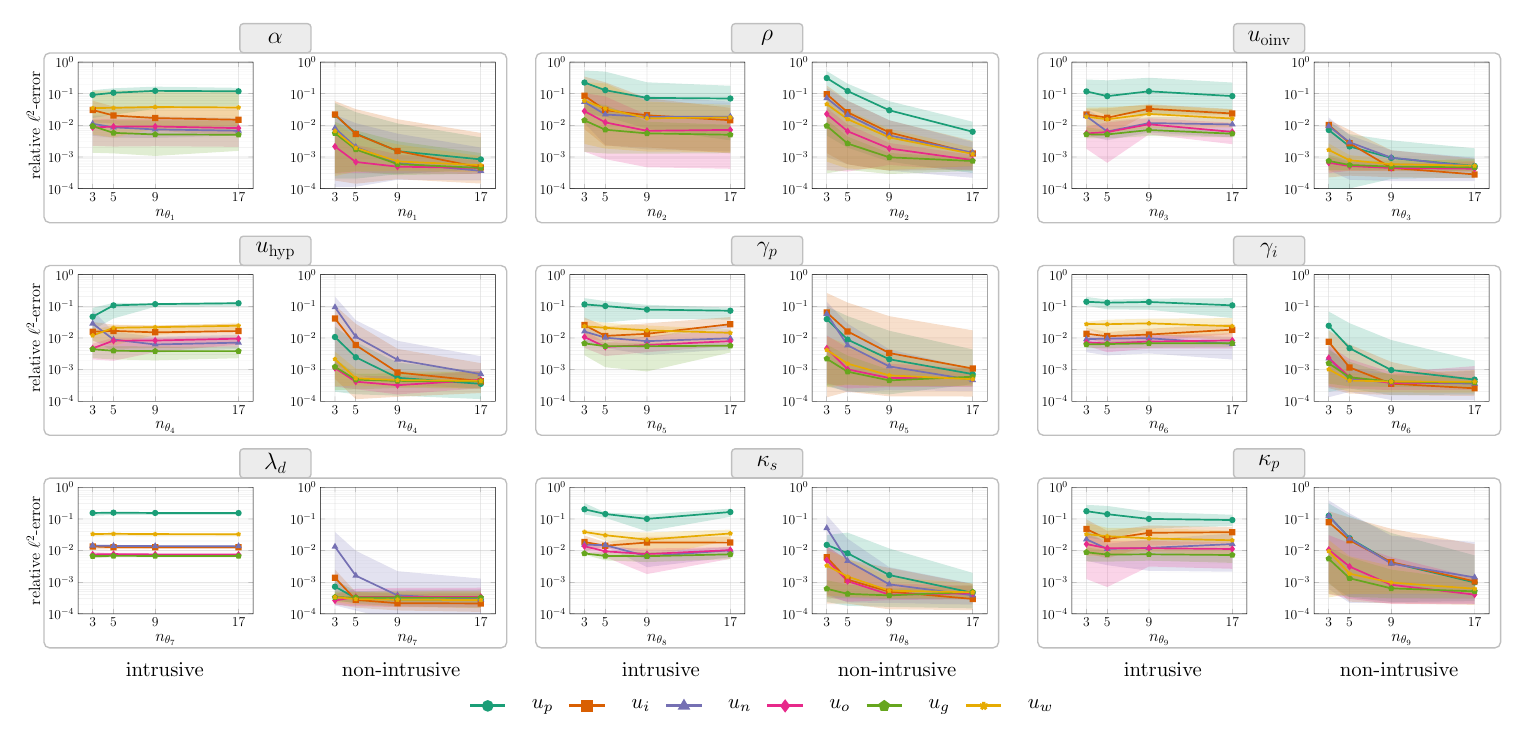}
\caption{Relative error with respect to the number of samples $n_{\theta_i}$ drawn to construct the  \htrom\ for each individual model parameter $\theta_i$. We consider the intrusive and non-intrusive variant of the \htrom. We report the relative  error compared to the \fom\ solution for an increasing number of samples  $n_{\theta_i}$ drawn. The samples are varied per parameter $\theta_i$, $i =  1,\ldots, 9$, with the  remaining parameters $\theta_j$, $j \not=i$ fixed. Each plot shows the trend  of the mean error averaged across 100 random off-grid samples $\theta_i$. The  envelopes correspond to the max and min error across all 100 trials. These results are for $n_t = 16$ using the tolerances reported in \Cref{t:rom-hyperparameters}.}
\label{f:error-htrom-parametersamples_nt16}
\end{figure}

\begin{figure}
\centering
\includegraphics[width=\textwidth]{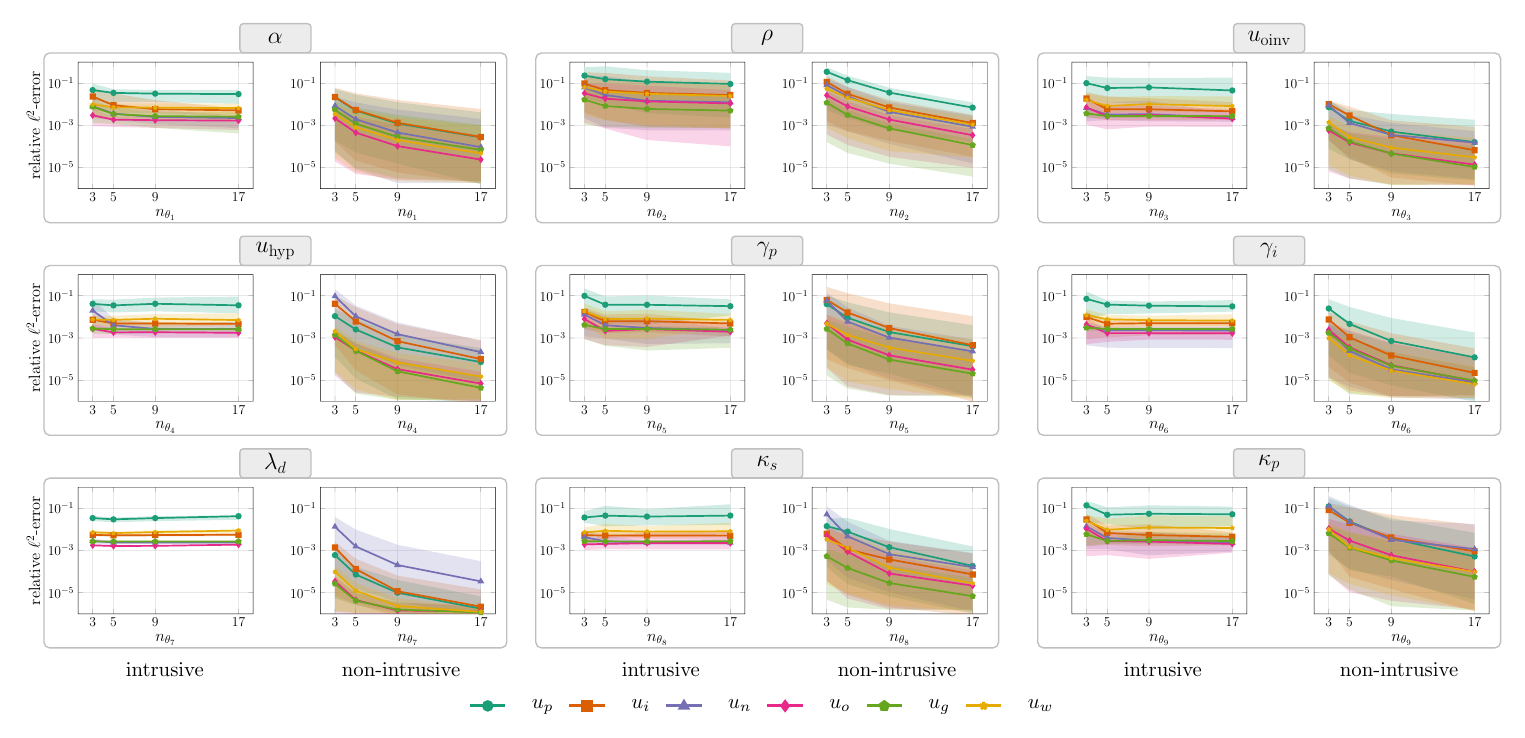}
\caption{Relative error with respect to the number of samples $n_{\theta_i}$ drawn to construct the  \htrom\ for each individual model parameter $\theta_i$. We consider the intrusive and non-intrusive variant of the \htrom. We report the error compared to the  \fom\ solution for an increasing number of samples $n_{\theta_i}$ drawn. The  samples are varied per parameter $\theta_i$. Each plot shows the trend of the  mean error averaged across 100 random off-grid samples $\theta_i$. The  envelopes correspond to the max and min error across all 100 trials. These  results are for $n_t = 64$ and $\epsilon^{\mathit{trom}}_{\mathit{on}}$, $\epsilon^{\mathit{trom}}_{\mathit{off}}$, and  $\epsilon^{\mathit{pre}}$ set to  1e--12.}
\label{f:error-htrom-parametersamples_nt64}
\end{figure}

\ipoint{Results} We report results for the intrusive and non-intrusive variants of the \htrom\ in \Cref{f:error-htrom-parametersamples_nt16,f:error-htrom-parametersamples_nt64}. We plot the mean error as a function of the number of samples $n_{\theta_i}$ for refining each individual parameter $\theta_i$ individually. Each plot shows the trend of the reconstruction error for each individual state variable. The envelopes show the maximum and minimum error obtained for each choice of $n_{\theta_i}$. \Cref{f:error-htrom-parametersamples_nt16} shows results for the parameters summarized in \Cref{t:rom-hyperparameters}. \Cref{f:error-htrom-parametersamples_nt64} includes results for the high-accuracy runs.

\ipoint{Observations} Overall, we have three sources for the reconstruction error for the intrusive model: a projection error (which should decrease as we increase the number of samples), an interpolation error (which should decrease as we increase the number of samples), and an error arising from the numerical time integration. Comparing the plots in \Cref{f:error-htrom-parametersamples_nt16,f:error-htrom-parametersamples_nt64}, we can observe that the error for the intrusive model seems to be dominated by the time integration error; the errors do not change significantly as we refine the mesh for the parameter samples. This is different for the results reported in \Cref{f:error-htrom-parametersamples_nt16}; the errors decrease as the mesh is refined, demonstrating that sampling at higher rates provides a smaller projection and/or interpolation error. Tightening the tolerances and increasing the number of time steps further reduces the errors (cf., \Cref{f:error-htrom-parametersamples_nt64}).

\subsection{Time Step Error}

\ipoint{Purpose} We investigate how the number of time steps $n_t$ used in the time integrator to compute the snapshot tensor $\mathcal{S}$ affects the \htrom\ accuracy.

\ipoint{Setup} We consider the number of time steps $n_t \in \{4,8,16\}$. This choice affects the accuracy of the snapshots stored in $\mathcal{S}$ during the offline phase of the \htrom.  For the intrusive variant, it also influences the accuracy in the online phase. To assess the impact, we compute the relative error between the \fom-solution and the \trom-solution at pseudo-time $t=1$. We select the number of time steps for the \fom-solver to $n_t = 128$ to obtain a high-accuracy reference solution. We randomly draw $10$ off-grid trial parameters $\theta^{\mathit{trial}} \in \Theta$ and report averages and maximum and minimum errors. We report errors for the intrusive and the non-intrusive variant of the \htrom.

\begin{figure}
\centering
\includegraphics[width=\textwidth]{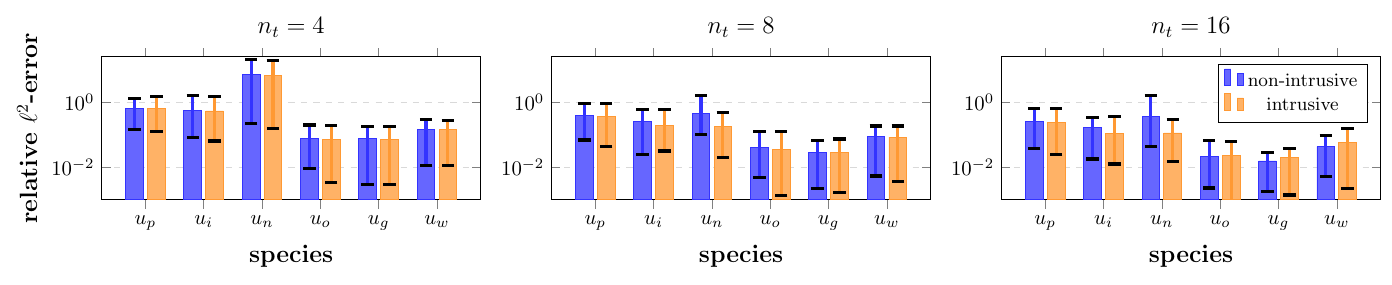}
\caption{Relative errors as a function of the number of time steps $n_t$ used in the \htrom\ with respect to \fom\ solution. We report the $\ell^2$-errors between a \fom\ high-fidelity reference solution for $n_t = 128$ and the \trom\ solution for a varying number of time steps $n_t$ used to compute the snapshots. We report errors with respect to the individual states $\up$, $\ui$, $\un$, $\uo$, $\ug$, and $\uw$. We report results for the intrusive and the non-intrusive \htrom\ variants. We run the \htrom\ for 10 randomly drawn off-grid model parameters $\theta^{\mathit{trial}} \in \Theta$. The bars represent the mean value across all 10 trials. The caps represent the min-max ranges of the error across all 10 trials.}
\label{f:error-num-timesteps}
\end{figure}

\ipoint{Results} We report the results in \Cref{f:error-num-timesteps}. We show the mean error as well as the min-max range for each state $\up$, $\ui$, $\un$, $\uo$, $\ug$, and $\uw$, respectively.

\ipoint{Observations} The most important observation is that the error is reduced as we increase the number of time steps. The error is most pronounced for the density $\un$ of necrotic tumor cells. There is no clear trend highlighting a difference between the non-intrusive (projection based) and intrusive \htrom\ variant. We select $n_t = 16$ as a default based on these experiments.

\subsection{Computational Performance as a Function of Mesh Refinement}

\ipoint{Purpose} We explore how mesh-refinement affects the computational performance of the \htrom.

\ipoint{Setup} We consider the multi-species model. We report runtime as a function of mesh resolution $(n_1,n_2,n_t)$ for the \htrom\ and the \fom-solver. We start with a mesh of resolution $(n_1,n_2,n_t) = (100,150,8)$ and increase the number of mesh points by a factor of two. To construct the \htrom\ we set the number of parameter samples to $n_{\theta_i} = 3$, $i=1,\ldots,9$. We construct the \htrom\ at the center point (on-grid location).

\ipoint{Results} We report the offline and online time for the \htrom\ and the runtime of the \fom-solver in \Cref{t:mesh-behavior}.

\begin{table}
\caption{Runtime (in seconds) for different stages of the proposed numerical framework as a function of the mesh size. We consider the multi-species model in two spatial dimensions. We report (from left to right) the time it takes to generate the snapshots, the remaining runtime of the offline stage (outside of snapshot generation), the evaluation of the \fom-solver for a trial parameter, and the online evaluation of the \htrom-solver variants (intrusive (IN) and non-intrusive (projection based; NI) variants).}
\label{t:mesh-behavior}
\tabadjust
\begin{tabular}{crrrrcc}
\toprule
\bf $n_1 \times n_2$ & $n_t$ & \bf snapshot & \bf offline & \bf \fom\ & \bf \htrom\ (IN) & \bf \htrom\ (NI) \\
\midrule
$100\times 150$ &  8 & \inum{03232.2} & \inum{ 814.0} & \fnum{0.13} & \fnum{0.09} & \snum{1.330e-3} \\
$200\times 300$ & 16 & \inum{23872.4} & \inum{1665.7} & \fnum{1.06} & \fnum{0.81} & \snum{7.239e-3} \\
\bottomrule
\end{tabular}
\end{table}

\ipoint{Observations} As we refine the mesh by a factor of 8 we can observe a close to linear increase in runtime for the individual stages of our numerical scheme.

\subsection{Comparison of {\rom}s}

Next, we compare the performance of the \fom-solver to the \podrom\ and different \trom\ variants.

\subsubsection{Two-Dimensional Multi-Species Model}

\ipoint{Purpose} We compare the performance of the \fom-solver to the \podrom\ and different \trom\ variants for the two-dimensional multi-species model.

\ipoint{Setup} We select the optimal algorithmic parameters as determined by the previous experiments. We consider the different \trom\ variants and the \podrom. We set the number of time steps for computing the snapshots to $n_t = 16$. To compare the models, we compute the optimal ranks for the \trom\ variants and the \podrom. We execute the \podrom\ using the optimal ranks as well as the optimal \trom\ ranks. The number of time steps for the \fom\ reference solution is also set to $n_t=16$.

We selected the tolerances to compute the online and offline ranks empirically.

To select the rank for the compression in space we explored tolerances between 1e--5 and 1e--8. We found that the reconstruction error remained at the same order for tolerances equal to or smaller than 1e--6. Consequently, we set $\epsilon_{\mathit{off}}^{\mathit{trom}}$ to 1e--6 for the space mode. We did similar experiments for the other modes. We set $\epsilon_{\mathit{off}}^{\mathit{trom}}$ to 1e--5 for the time mode and 1e--7 for the parameter modes. We summarize these choices in \Cref{t:rom-hyperparameters}.

\begin{figure}
\centering
\includegraphics[width=0.8\textwidth]{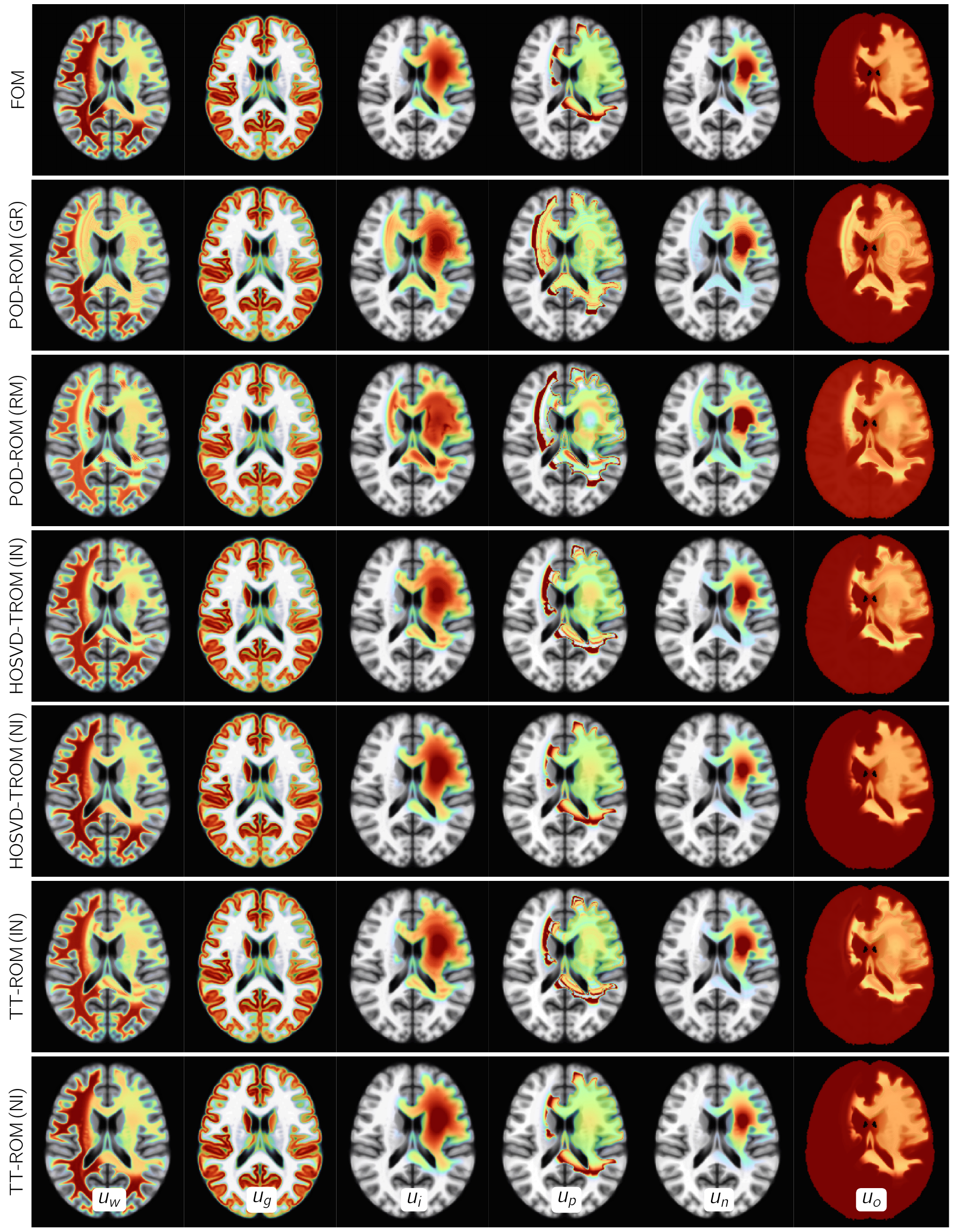}
\caption{We show simulation results comparing the \fom\ to different \rom\  variants. These results are for the multispecies model (two-dimensional  simulations). We show (from left to right) axial slices for the white matter  density $\uw$, the gray matter density $\ug$, and the densities for the  infiltrative tumor cells $\ui$, proliferative tumor cells $\up$, necrotic  tumor cells $\un$, and the oxygen concentration $\uo$. The top row shows  solutions for the \fom. The second row shows solutions for the \podrom\  (global rank selection). The determined ranks are 140, 45, 132, 260, 140, and  35 for $\uw$, $\ug$, $\ui$, $\up$, $\un$ and $\uo$, respectively. The third  row shows solutions for the \podrom, where we matched the \podrom\ ranks with  the online ranks for the \htrom. The individual ranks are 10, 6, 12, 16, 11,  and 8 for $\uw$, $\ug$, $\ui$, $\up$, $\un$ and $\uo$, respectively. The  fourth and fifth rows show results for the \htrom\ intrusive (IN) and  non-intrusive (NI) variants. The last two rows shows the solution of \ttrom\  for the intrusive (IN) and non-intrusive (NI) variants with ranks 10, 6, 12,  16, 11 and 8 for $\uw$, $\ug$, $\ui$, $\up$, $\un$ and $\uo$, respectively. We  show the associated local error maps in  \Cref{f:rom-vs-fom-multispecies-2D-errormaps}. The online and offline ranks are reported in \Cref{t:optimal-ranks-hosvd-rom,t:optimal-ranks-tt-rom}, respectively.}
\label{f:rom-vs-fom-multispecies-2D-solutions}
\end{figure}

\begin{figure}
\centering
\includegraphics[width=0.8\textwidth]{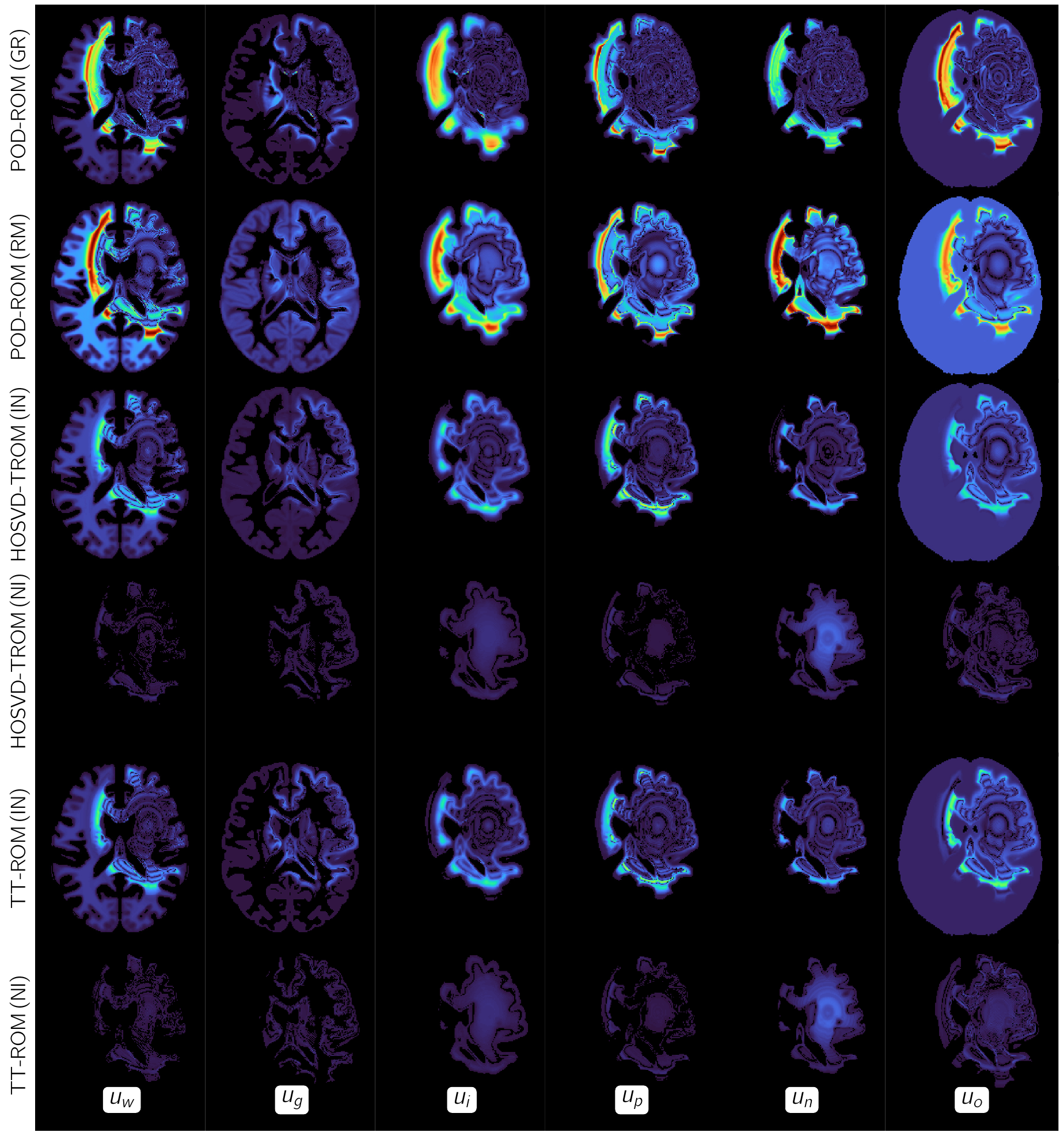}
\caption{Error maps for the different \rom\ variants. We compute the error  between the considered \rom\ variants and a high-fidelity \fom\ solution.  These results are for the multispecies model (two-dimensional simulations). We  show (from left to right) axial slices for the white matter density $\uw$, the  gray matter density $\ug$, and the densities for the infiltrative tumor cells  $\ui$, proliferative tumor cells $\up$, necrotic tumor cells $\un$, and the  oxygen concentration $\uo$. To be able to compare the error maps for the  various approaches, we normalized the color maps to the same range. The first  and second rows show pointwise error maps for the \podrom\ for both global  rank selection and rank matched to the \htrom, respectively. The middle two  rows show pointwise error maps for the \htrom\ for intrusive (IN) and  non-intrusive (NI) variants. The bottom two rows show pointwise error maps for  the \ttrom\ for both intrusive (IN) and non-intrusive (NI) variants,  respectively. The results correspond to those reported in  \Cref{f:rom-vs-fom-multispecies-2D-solutions}. The online and offline ranks  are reported in \Cref{t:optimal-ranks-hosvd-rom,t:optimal-ranks-tt-rom}, respectively.}
\label{f:rom-vs-fom-multispecies-2D-errormaps}
\end{figure}

\ipoint{Results} The optimal ranks for the \htrom\ are reported in \Cref{t:optimal-ranks-hosvd-rom}. The optimal ranks for the \ttrom\ are reported in \Cref{t:optimal-ranks-tt-rom}. We show exemplary solution in \Cref{f:rom-vs-fom-multispecies-2D-solutions}. We show the associated pointwise error maps in \Cref{f:rom-vs-fom-multispecies-2D-errormaps}. We report quantitative results in \Cref{t:relative-ell2-error-roms-vs-fom-ms2d} (relative $\ell^2$-error). Runtimes for the online stage are reported in \Cref{t:runtimes-roms-vs-fom-ms2d}.

\begin{table}
\caption{Optimal ranks for the \htrom. The results are for the two-dimensional multi-species model. We report the ranks determined with respect to each state variable. The associated ranks (in the same order of species) for the \podrom\ are 140, 45, 132, 260, 140, and 35, respectively.}
\label{t:optimal-ranks-hosvd-rom}
\tabadjust
\begin{tabular}{lrrrrrrrrrrrr}
\toprule
\bf state & \multicolumn{11}{c}{\bf offline ranks} & \bf online ranks\\
\midrule
$\uw$ & 140 & 10 & 3 & 10 & 3 & 3 & 3 & 5 & 5 & 3 & 3 & 10  \\
$\ug$ &  45 &  6 & 3 &  9 & 3 & 3 & 3 & 5 & 4 & 3 & 3 &  6  \\
$\ui$ & 132 & 12 & 3 & 10 & 3 & 3 & 3 & 5 & 5 & 3 & 3 & 12  \\
$\up$ & 260 & 16 & 3 & 10 & 3 & 3 & 3 & 5 & 5 & 3 & 3 & 16  \\
$\un$ & 140 & 14 & 3 & 10 & 3 & 3 & 3 & 5 & 5 & 3 & 3 & 11  \\
$\uo$ &  35 &  8 & 3 &  9 & 3 & 3 & 3 & 5 & 4 & 3 & 3 &  8  \\
\bottomrule
\end{tabular}
\end{table}

\begin{table}
\caption{Optimal ranks for the \ttrom. The results are for the two-dimensional multi-species model. We report the ranks determined with respect to each state variable. The associated ranks (in the same order of species) for the \podrom\ are 140, 45, 132, 260, 140, and 35, respectively.}
\label{t:optimal-ranks-tt-rom}
\tabadjust
\begin{tabular}{lrrrrrrrrrrrrr}
\toprule
\bf state & \multicolumn{12}{c}{\bf offline ranks} & \bf online ranks\\\midrule
$\uw$ & 1 & 140 &  74 & 120 & 273 & 285 & 211 & 119 & 35 & 9 & 3 & 1 & 10 \\
$\ug$ & 1 &  45 &  14 &  22 &  22 &  25 &  20 &  20 & 13 & 7 & 3 & 1 &  6 \\
$\ui$ & 1 & 132 & 133 & 235 & 321 & 318 & 246 & 163 & 37 & 9 & 3 & 1 & 12 \\
$\up$ & 1 & 260 & 308 & 655 & 775 & 604 & 346 & 171 & 40 & 9 & 3 & 1 & 16 \\
$\un$ & 1 & 140 & 183 & 363 & 526 & 444 & 278 & 156 & 39 & 9 & 3 & 1 & 11 \\
$\uo$ & 1 &  35 &  16 &  23 &  26 &  34 &  29 &  29 & 17 & 8 & 3 & 1 &  8 \\
\bottomrule
\end{tabular}
\end{table}

\begin{table}
\caption{Relative $\ell^2$-errors at final time $t=1$. We compare different  \rom\ variants to the \fom. These results are for the two-dimensional  multi-species tumor growth model. We report (from left to right) errors for  the \podrom\ with optimal rank (GR), the \podrom\ with rank matched to the  \htrom\ (RM), and errors for the intrusive \htrom\ variant (IN), the  non-intrusive \htrom\ variant (NI), the intrusive \ttrom\ variant (IN), and  the non-intrusive \ttrom\ variant (NI), respectively. The online and offline  ranks are reported in \Cref{t:optimal-ranks-hosvd-rom,t:optimal-ranks-tt-rom}.}
\label{t:relative-ell2-error-roms-vs-fom-ms2d}
\tabadjust
\begin{tabular}{ccccccc}
\toprule
& \multicolumn{2}{c}{\bf \podrom} & \multicolumn{2}{c}{\bf \htrom} & \multicolumn{2}{c}{\bf \ttrom} \\
\bf state & \bf GR & \bf RM & \bf IN & \bf NI & \bf IN & \bf NI \\
\midrule
$\uw$ & \snum{1.843e-01} & \snum{2.510e-01} & \snum{9.824e-02} & \snum{5.002e-03} & \snum{8.066e-02} & \snum{6.068e-03}\\
$\ug$ & \snum{2.647e-02} & \snum{3.748e-02} & \snum{1.908e-02} & \snum{3.315e-03} & \snum{2.054e-02} & \snum{4.062e-03}\\
$\ui$ & \snum{4.722e-01} & \snum{5.104e-01} & \snum{1.652e-01} & \snum{2.909e-02} & \snum{1.501e-01} & \snum{2.835e-02}\\
$\up$ & \snum{8.753e-01} & \snum{1.002e+00} & \snum{5.619e-01} & \snum{2.850e-02} & \snum{4.413e-01} & \snum{3.173e-02}\\
$\un$ & \snum{2.557e-01} & \snum{5.100e-01} & \snum{9.096e-02} & \snum{6.445e-02} & \snum{8.479e-02} & \snum{6.137e-02}\\
$\uo$ & \snum{9.215e-02} & \snum{9.036e-02} & \snum{3.750e-02} & \snum{3.639e-03} & \snum{3.806e-02} & \snum{6.007e-03}\\
\bottomrule
\end{tabular}
\end{table}

\begin{table}
\caption{Runtime comparison (online stage). We report results for the  two-dimensional multi-species tumor growth model. We report runtimes (from top  to bottom) for the \fom, the \podrom\ with optimal rank (GR), the \podrom\  with rank matched to the \htrom\ (RM), the intrusive \htrom\ variant (IN), the  non-intrusive \htrom\ variant (NI), the intrusive \ttrom\ variant (IN), and  non-intrusive \ttrom\ variant (NI). The reported runtimes are in seconds and  for the online phase. The ranks for constructing the {\rom}s are reported in \Cref{t:optimal-ranks-hosvd-rom,t:optimal-ranks-tt-rom}.}
\label{t:runtimes-roms-vs-fom-ms2d}
\tabadjust
\begin{tabular}{lrr}
\toprule
\bf approach & \bf runtime & \bf speedup \\
\midrule
FOM           & \fnum{1.110} & --- \\
\podrom\ (GR) & \fnum{6.732} & \fnum{ 0.165} \\
\podrom\ (RM) & \fnum{0.713} & \fnum{ 1.556} \\
\htrom\ (IN)  & \fnum{0.795} & \fnum{ 1.400} \\
\htrom\ (NI)  & \fnum{0.013} & \fnum{85.380} \\
\ttrom\ (IN)  & \fnum{0.851} & \fnum{ 1.300} \\
\ttrom\ (NI)  & \fnum{0.013} & \fnum{ 85.38} \\
\bottomrule
\end{tabular}
\end{table}

\ipoint{Observations} The most important observation is that the non-intrusive (projection based) \trom\ variants outperform the intrusive variants, with a speedup of up to 85$\times$ and an excellent agreement with the \fom-solution across all species. The \htrom\ and the \ttrom\ display a very similar performance in terms of accuracy (see \Cref{t:relative-ell2-error-roms-vs-fom-ms2d}) and runtime (see \Cref{t:runtimes-roms-vs-fom-ms2d}).

This is expected, since we select the same tolerances for the construction of the \htrom\ and \ttrom. We note that while the ranks for these \trom\ representations are different, with higher offline ranks observed for the \ttrom\ (see \Cref{t:optimal-ranks-tt-rom}), the ranks of these two variants have a quite different meaning. In general, the complexity of the \ttrom\ representation is much smaller; the complexity of the \htrom\ is multiplicative in the ranks whereas the complexity of the \ttrom\ is additive across the individual cores.

The runtime performance of the \podrom\ (rank matched) and the intrusive \trom\ variants are similar. To achieve a good accuracy for the \podrom\ we had to use a small tolerance; the compression for the optimal rank was not significant enough to obtain a competitive runtime, even compared to the \fom-solver. If we relax the energy tolerance for the \podrom, the runtime is faster than the runtime for the \fom-solver. We also note that we did not observe this deterioration in runtime performance for the three-dimensional case reported in the next section.

\subsubsection{Three-Dimensional Single-Species Model}

\ipoint{Purpose} We compare the performance of the \fom-solver and different \rom\ variants for the 3D implementation of the single species model.

\ipoint{Setup} We consider the \podrom\ and the \htrom. We determined the optimal ranks for the {\rom}s. We execute the time integrator for $n_t = 16$.

\begin{figure}
\centering
\includegraphics[width=0.8\textwidth]{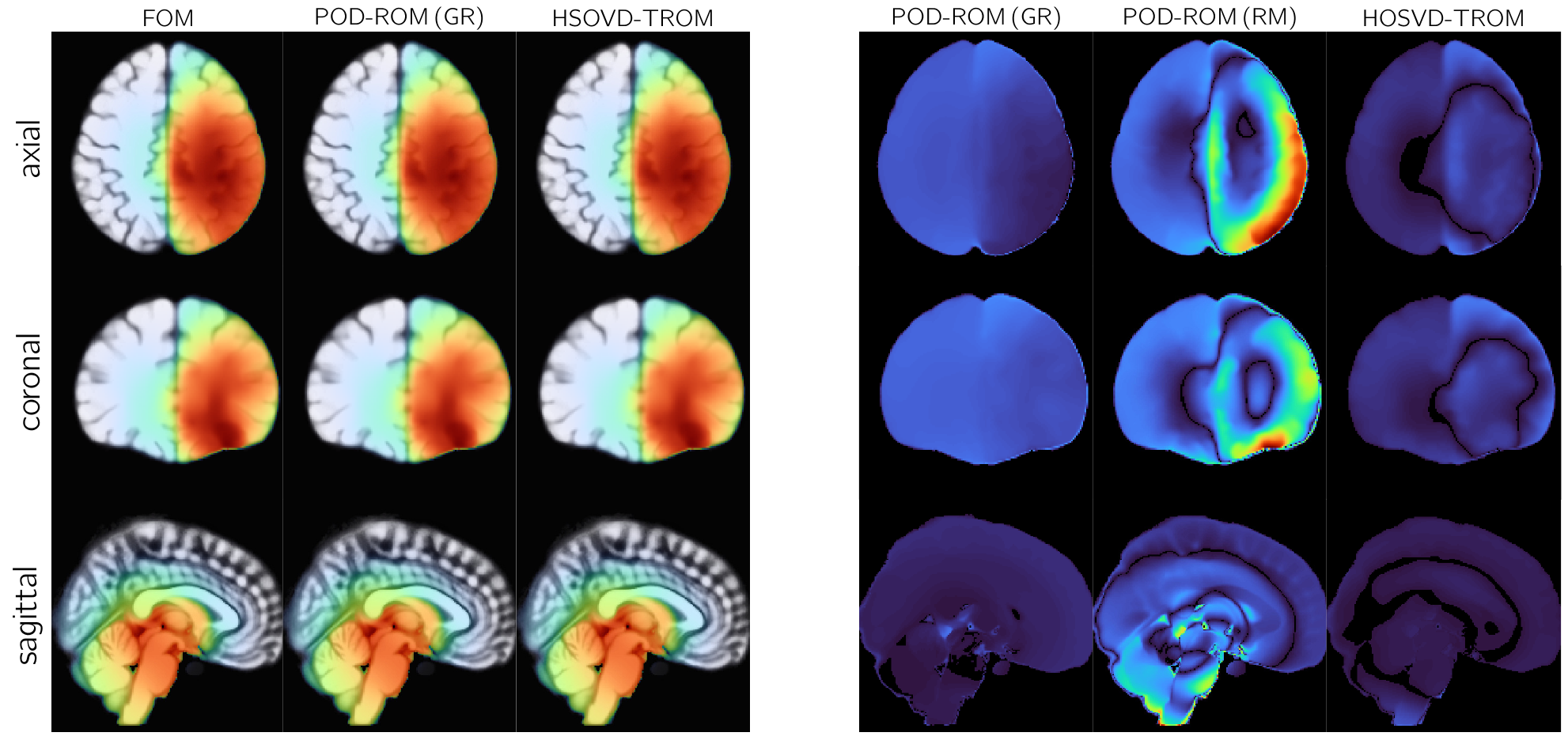}
\caption{Visualization of simulation results for the 3D single species model. We show (from top to bottom) axial, coronal, and sagittal views. The left block shows simulation results. The right block shows error maps (normalized to the same range). Each column corresponds to a different solver; from left to right: \fom\ solution, \podrom\ solution with optimal rank (GR) and rank matched to the \htrom\ rank (RM), and the \htrom\ (non-intrusive variant). We report quantitative results in \Cref{t:results-3d-single-species}.}
\label{f:results-3d-single-species}
\end{figure}

\ipoint{Results} The ranks for the \podrom\ and the \htrom\ are reported in \Cref{t:ranks-rom-pod-v-trom-ss3d}. We showcase simulation results in \Cref{f:results-3d-single-species}. The reconstruction error, runtimes, and speedups are reported in \Cref{t:results-3d-single-species}.

\begin{table}
\caption{Optimal ranks for the \podrom\ and the \htrom\ for the three-dimensional single species model. For the \htrom\ we report the estimated ranks for space, for time, and the diffusion coefficient $\alpha$.}
\label{t:ranks-rom-pod-v-trom-ss3d}
\tabadjust
\begin{tabular}{ccc}
\toprule
\bf \podrom\ rank & \multicolumn{2}{c}{\bf \htrom} \\
& \bf offline ranks & \bf online rank \\
\midrule
48  & 48 \,\,\, 14 \,\,\, 5 & 10 \\
\bottomrule
\end{tabular}
\end{table}

\begin{table}
\caption{Relative error at final time $t=1$ and runtime for the 3D single-species tumor growth model. We report results for the \podrom\ (optimal rank; GR), \podrom\ (rank matched with the \htrom; RM), and the \htrom\ (for the intrusive (IN) and non-intrusive (NI) variants). We report the relative $\ell^2$-error and the $\ell^\infty$-error between the \rom\ solution and the \fom\ solution. We also report the runtime (in seconds) and the speedup compared to the \fom\ solver. We show simulation results in \Cref{f:results-3d-single-species}.}
\label{t:results-3d-single-species}
\tabadjust
\begin{tabular}{lrrrr}
\toprule
\bf model & \bf $\ell^2_{\text{rel}}$-error & \bf $\ell^{\infty}$-error  & \bf runtime & \bf speedup \\
\midrule
\fom          & ---               & ---              & \fnum{345.230} & ---  \\
\podrom\,(GR) & \snum{1.086e-03}  & \snum{3.031e-03} & \fnum{ 42.683} & \fnum{   8.088}$\times$ \\
\podrom\,(RM) & \snum{7.577e-03}  & \snum{2.487e-02} & \fnum{ 13.586} & \fnum{ 25.4107}$\times$ \\
\htrom\,(IN)  & \snum{7.362e-04}  & \snum{3.432e-03} & \fnum{ 14.124} & \fnum{ 24.4427}$\times$ \\
\htrom\,(NI)  & \snum{7.826e-04}  & \snum{2.639e-03} & \fnum{  2.833} & \fnum{121.8602}$\times$ \\
\bottomrule
\end{tabular}
\end{table}

\ipoint{Observations} The most important observation is that the {\rom}s deliver results that agree with the \fom\ solution to high accuracy with speedups ranging from 8$\times$ for the \podrom\ up to 120$\times$ for the projection based \htrom.

\section{Conclusions}\label{s:conclusions}

We have designed a numerical framework that takes advantage of {\trom}s for the efficient simulation of brain tumor growth. Our framework integrates a finite volume discretization with an operator-splitting time integration strategy and a family of projection-based surrogate models---including a \podrom\ baseline and two tensorial variants (\htrom\ and \ttrom), each available in intrusive and non-intrusive implementations. A spatial pre-compression stage further reduces the memory and computational footprint of the offline phase, making the approach tractable for large three-dimensional problems. The most important observations are as follows.
\begin{itemize}[leftmargin=*]
\item The {\trom}s outperform the \fom-solver and the {\podrom} variant, delivering excellent accuracy with speedups ranging from $85\times$ to $120\times$ for the two- and three-dimensional models considered in the present study.
\item Our framework allows us to effectively handle complex, coupled PDE systems controlled by up to $m=9$ model parameters.
\item The operator-splitting strategy allows us to efficiently deploy {\rom}s by separating linear and non-linear operators, allowing for a straightforward application of {\rom} approaches.
\item The \htrom\ and \ttrom\ surrogates both outperform the classical \podrom\ baseline. The tensorial structure enables simultaneous compression across the spatial, temporal, and parameter dimensions. Both intrusive and non-intrusive variants deliver competitive accuracy; the non-intrusive variant offers a particularly lightweight online stage that avoids reduced-order time integration entirely.
\item The spatial pre-compression stage is essential for tractability in three-dimensional settings: it reduces the leading spatial dimension of the snapshot tensor from $n_x$ to $r_x \ll n_x$ prior to the tensor decomposition, enabling offline construction on standard computing hardware.
\item The proposed framework accommodates both the single-species reaction--diffusion formulation and the more complex six-state, nine-parameter go-or-grow multi-species model, demonstrating its applicability to parametric PDE systems of varying complexity.
\end{itemize}

Taken together, these results demonstrate that projection-based tensorial surrogates are a compelling alternative to purely data-driven approaches for parametric biomedical simulations. They provide rigorous approximation guarantees, preserve the physical structure of the governing equations, and achieve substantial online speedups at a fraction of the cost of the \fom. The offline investment---repeated high-fidelity solves over the parameter training set---is amortized over the many online evaluations required in downstream many-query settings.

Our future work will target the deployment of our solver to dedicated hardware architectures and more effective programming languages and models, with the ultimate aspiration of integrating the proposed framework into many-query settings. Specifically, we aim to embed the proposed surrogates within Bayesian inference and inverse problem solvers for patient-specific parameter estimation from medical imaging data, and within digital-twin workflows for treatment planning and outcome prediction. Further directions include extending the framework to higher-dimensional parameter spaces through adaptive sampling strategies and exploring the combination of tensorial surrogates with data-driven corrections to improve predictive accuracy in extrapolation regimes.

\paragraph{Acknowledgements} This work was partly supported by the National Science Foundation ({\bf NSF}) under the award DMS-2145845 (AI \& AM) and  DMS-2309197 (RM \& MO). Any opinions, findings, and conclusions or recommendations expressed herein are those of the authors and do not necessarily reflect the views of NSF.  This work was completed while AM was in residence at the Institute for Computational and Experimental Research in Mathematics in Providence, RI, during the ``Stochastic and Randomized Algorithms in Scientific Computing: Foundations and Applications'' program.  We are grateful for the support of the Research Computing Data Core at the University of Houston.

\begin{appendix}

\section{Hardware}

All runs for snapshot generation and \rom\ training/evaluation were executed on the \emph{Sabine} cluster operated by the University of Houston Research Computing Data Core, using the SLURM workload manager. The jobs were run as single-node CPU jobs with the following resource requests (nodes: 1; tasks: 1; CPU cores per task: 32; memory: 64 GB).

\section{Software and Libraries}

Code development was done using MATLAB \texttt{R2025a}. The implementation uses sparse linear algebra primitives as well as low-level parallel computing primitives for distributed snapshot generation. \htrom\ and \ttrom\ operations in the offline stage use the MATLAB Tensor Toolbox (\texttt{tensor}, \texttt{ttm})~\cite{bader2008:efficient, kolda2017:tensor, tensortoolbox2025}.

\section{Data}

For the brain geometry we consider the MNI ICBM152 templates, in particular, the
  ``ICBM 2009c Nonlinear Symmetric'' template~\cite{fonov2011:unbiased, fonov2009:unbiased, collins1999:animal}. The image is a three-dimensional volume which provides probability maps for different tissue compartments as well as anatomical images.

\section{Algorithms}

  We outline the main algorithmic approaches proposed in the present work.  \Cref{a:rsvd-fixed-accuracy} illustrates the considered implementation of the  rSVD algorithm using an adaptive rank finder strategy. The different stages and variants of the \htrom\ are presented in \Cref{a:hosvd-trom-offline,a:hosvd-trom-online,a:hosvd-trom-online-ni}. The different stages and variants of the \ttrom\ are presented in \Cref{a:tt-trom-offline,a:tt-trom-online,a:tt-trom-online-ni}. The offline precompression is outlined in \Cref{a:precompression-offline}.

\begin{algorithm}
\caption{Randomized SVD ({\bf rSVD}) using an adaptive rank finder strategy~\cite{halko2011:finding}.}
\label{a:rsvd-fixed-accuracy}
\begin{algorithmic}[1]
\State \textbf{input:} $A \in \mathbb{R}^{\tilde{n}\times \tilde{m}}$
\State \textbf{parameters:} tolerance $\epsilon>0$; oversampling $p_o \in \mathbb{N}_0$; target rank $r \in \mathbb{N}$; growth step $p_s \in \mathbb{N}$; power iterations $q \in \mathbb{N}_0$
\While{true}
    \State $\ell \gets \min(r + p_o, \min(\tilde{n},\tilde{m}))$\label{l:oversampling}
    \State draw $\Omega \sim \mathcal{N}(0,1)^{\tilde{m}\times \ell}$
    \State form sketch $Y \gets A\Omega$
    \For{$s=1, \ldots, q$}
        \State $Q \gets \operatorname{orth}(Y)$
        \State $Y \gets A(A^\mathsf{T}Q)$
    \EndFor
    \State $Q \gets \operatorname{orth}(Y)$
    \State $B \gets Q^\mathsf{T} A$
    \State $(\tilde{U},\Sigma,V) \gets \operatorname{svd}(B)$
    \State compute cumulative energy $\eta(r) \gets \sum_{i=1}^r (\Sigma_{ii})^2 \big/ \|A\|_F^2$ (cf. \Cref{e:energy-bound-r})
    \If{$\eta(r) \ge 1-\epsilon$}
        \State \Return $U = Q\,\tilde{U}(\,:\,,1:r)$, $\Sigma=\Sigma(1:r,1:r)$, $V=V(\,:\,,1:r)$
    \ElsIf{$\ell = \min(N,M)$}
        \State \Return $U = Q\,\widetilde U$, $\Sigma$, $V$
    \Else
        \State $r \gets \min(r + p_s, \min(\tilde{n},\tilde{m}))$\label{l:range-finder} 
    \EndIf
\EndWhile
\State \textbf{output:} $U \in \mathbb{R}^{\tilde{n} \times r}$, $\Sigma \in \mathbb{R}^{r\times r}$, $V \in \mathbb{R}^{\tilde{m}\times r}$
\end{algorithmic}
\end{algorithm}

\begin{algorithm}
\caption{Offline stage of the interpolatory \htrom.\label{a:hosvd-trom-offline}}
\begin{algorithmic}[1]
\State \textbf{input:} tolerance $\epsilon^{\mathit{trom}}_{\mathit{off}} \in (0,1)$; parameter samples $\Theta_s$
\State compute snapshot tensor $\mathcal{S}$ by evaluating \fom\ for each $\theta^l \in \Theta_s$, $l \in \mathcal{I}_s$
\For{$i=1, \ldots, m+2$}
    \State $S_{(i)} \gets$ perform mode-$i$ unfolding of tensor $\mathcal{S}$
    \State compute $S_{(i)} = U^{(i)}_{r_i} \Sigma^{(i)}_{r_i}(V^{(i)}_{r_i})^{\mathsf{T}}$ with $r_i(\epsilon^{\mathit{trom}}_{\mathit{off}})$ chosen via \Cref{e:energy-bound-r}
    \State $\Phi_{(i)} \gets U^{(i)}_{r_i}$
\EndFor
\State $\Phi^{(x)} \gets \Phi_{(1)}$, $\Phi^{(t)} \gets \Phi_{(2)}$, $\Phi^{(\theta_i)} \gets \Phi_{(i+2)}$, $i=1,\dots,m$
\State compute core tensor $\mathcal{C}$ as in \Cref{e:core-tensor}
\State \textbf{output:} $\{\mathcal{C}, \Phi^{(x)}, \Phi^{(t)}, \Phi^{(\theta_1)}, \ldots, \Phi^{(\theta_m)}\}$
\end{algorithmic}
\end{algorithm}

\begin{algorithm}
\caption{Online stage for the intrusive variant of the interpolatory \htrom.\label{a:hosvd-trom-online}}
\begin{algorithmic}[1]
\State \textbf{input:} $\{\mathcal{C}, \Phi^{(x)}, \Phi^{(t)}, \Phi^{(\theta_1)}, \ldots, \Phi^{(\theta_m)}\}$; tolerance $\epsilon^{\mathit{trom}}_{\mathit{on}} \in (0,1)$, $\theta \in \mathbb{R}^m$
\State assemble core matrix $\tilde{C}(\theta)$ based on core tensor contraction \Cref{e:core-tensor-contraction}
\State $[U_r\Sigma_rV_r^\mathsf{T}](\theta) \gets$ compute truncated SVD of $\tilde{C}(\theta)$ with $r(\epsilon^{\mathit{trom}}_{\mathit{on}})$ chosen based on \Cref{e:energy-bound-r}
\State $\Psi(\theta) \gets \Phi^{(x)}(\theta)U_r(\theta)$
\State $\mu_h \gets$ numerically integrate \Cref{e:reduced-system-trom}
\State $u_h^j(\theta) \gets \Psi(\theta) \mu_h^j(\theta)$ for all $j = 1,\ldots,n_t$
\State \textbf{output:} $u_h^j(\theta)$
\end{algorithmic}
\end{algorithm}

\begin{algorithm}
\caption{Online stage of the non-intrusive variant of the interpolatory \htrom.\label{a:hosvd-trom-online-ni}}
\begin{algorithmic}[1]
\State \textbf{input:} $\{\mathcal{C}, \Phi^{(x)}, \Phi^{(t)}, \Phi^{(\theta_1)}, \ldots, \Phi^{(\theta_m)}\}$; tolerance $\epsilon^{\mathit{trom}}_{\mathit{on}} \in (0,1)$, $\theta \in \mathbb{R}^m$
\State assemble core matrix $\tilde{C}(\theta)$ based on core tensor contraction \Cref{e:core-tensor-contraction}
\State $[U_r\Sigma_rV_r^\mathsf{T}](\theta) \gets$ compute truncated SVD of $\tilde{C}(\theta)$ with $r(\epsilon^{\mathit{trom}}_{\mathit{on}})$ chosen based on \Cref{e:energy-bound-r}
\State $\Psi^{(x)}(\theta) \gets \Phi^{(x)} U_r(\theta)$
\State $\Psi^{(t)}(\theta) \gets \Phi^{(t)} V_r(\theta)$
\State $\hat{X}(\theta) \gets \Psi^{(x)}(\theta)\,\Sigma_r(\theta) \bigl(\Psi^{(t)}(\theta)\bigr)^\mathsf{T}$
\State $u_h^j(\theta) \gets \hat{X}(\theta)_{:,j}$ for all time indices $j = 1,\ldots,n_t$
\State \textbf{output:} $u_h^j(\theta)$
\end{algorithmic}
\end{algorithm}

\begin{algorithm}
\caption{Offline stage of interpolatory \ttrom\ (TT-SVD construction).\label{a:tt-trom-offline}}
\begin{algorithmic}[1]
\State \textbf{input:} tolerance $\epsilon^{\mathit{trom}}_{\mathit{off}} \in (0,1)$ and parameter samples $\Theta_s$
\State compute snapshot tensor $\mathcal{S}$ by evaluating \fom\ for each $\theta^l\in\Theta_s$, $l\in\mathcal{I}_s$
\State set $n_1\gets n_x$, $n_2\gets n_t$, and $n_{k+2}\gets n_{\theta_k}$ for $k=1,\ldots,m$; initialize $r_0\gets 1$ and $\mathcal{Q}^{(1)}\gets\mathcal{S}$
\For{$k = 1,\ldots,m+1$}
    \State view $\mathcal{Q}^{(k)}$ as an element of $\mathbb{R}^{r_{k-1}\times n_k\times n_{k+1}\times\cdots\times n_{m+2}}$
    \State form the compound unfolding $Q^{(k)}\in\mathbb{R}^{(r_{k-1}n_k)\times(n_{k+1}\cdots n_{m+2})}$
    \State $U^{(k)}_{r_k}\,\Sigma^{(k)}_{r_k}\,\bigl(V^{(k)}_{r_k}\bigr)^\mathsf{T} \gets$ compute thin SVD of $Q^{(k)}$ with $r_k(\epsilon^{\mathit{trom}}_{\mathit{off}})$ chosen via \Cref{e:energy-bound-r}
    \State define the $k$th TT core $\mathcal{G}^{(k)}\in\mathbb{R}^{r_{k-1}\times n_k\times r_k}$ by reshaping $U^{(k)}_{r_k}$ as in \Cref{e:tt-core-from-U}
    \State $\mathcal{Q}^{(k+1)} \gets \texttt{reshape}\bigl(\Sigma^{(k)}_{r_k}\,\bigl(V^{(k)}_{r_k}\bigr)^\mathsf{T},\,(r_k, n_{k+1}, \ldots, n_{m+2})\bigr)$
\EndFor
\State define the final TT core $\mathcal{G}^{(m+2)}\in\mathbb{R}^{r_{m+1}\times n_{m+2}\times 1}$ by reshaping $\mathcal{Q}^{(m+2)}$ as in \Cref{e:tt-final-core}
\State \textbf{output:} $\{\mathcal{G}^{(k)}\}_{k=1}^{m+2}$
\end{algorithmic}
\end{algorithm}

\begin{algorithm}
\caption{Online stage of intrusive variant of interpolatory \ttrom.\label{a:tt-trom-online}}
\begin{algorithmic}[1]
\State \textbf{input:} TT cores $\{\mathcal{G}^{(k)}\}_{k=1}^{m+2}$; tolerance $\epsilon^{\mathit{trom}}_{\mathit{on}}\in(0,1)$; $\theta \in \mathbb{R}^m$
\State $\Phi^{(x)} \gets \texttt{reshape}(\mathcal{G}^{(1)}, (n_x, r_1))$ (mode-$2$ unfolding of first TT core; $r_0=1$)
\State assemble contracted matrix $\tilde{C}(\theta')$ based on TT contraction \Cref{e:tt-param-core-contraction}--\Cref{e:tt-core-tensor-contraction}
\State $[U_r\Sigma_rV_r^\mathsf{T}](\theta) \gets$ compute thin SVD of $\tilde{C}(\theta)$ with $r(\epsilon^{\mathit{trom}}_{\mathit{on}})$ chosen via \Cref{e:energy-bound-r}
\State $\Psi(\theta) \gets \Phi^{(x)}(\theta) U_r(\theta)$
\State $\mu_h \gets$ numerically integrate \Cref{e:reduced-system-trom}
\State $u_h^j(\theta)\gets \Psi(\theta)\mu_h^j(\theta)$ for all $j = 1,\ldots,n_t$
\State \textbf{output:} $u_h^j(\theta)$
\end{algorithmic}
\end{algorithm}

\begin{algorithm}
\caption{Online stage of non-intrusive variant of interpolatory \ttrom.\label{a:tt-trom-online-ni}}
\begin{algorithmic}[1]
\State \textbf{input:} TT cores $\{\mathcal{G}^{(k)}\}_{k=1}^{m+2}$; tolerance $\epsilon^{\mathit{trom}}_{\mathit{on}}\in(0,1)$; $\theta \in \mathbb{R}^m$
\State $\Phi^{(x)} \gets \texttt{reshape}(\mathcal{G}^{(1)}, (n_x, r_1))$ (mode-$2$ unfolding of first TT core; $r_0=1$)
\State assemble contracted matrix $\tilde{C}(\theta)$ based on TT contraction \Cref{e:tt-param-core-contraction}--\Cref{e:tt-core-tensor-contraction}
\State $[U_r\Sigma_rV_r^\mathsf{T}](\theta) \gets$ compute thin SVD of $\tilde{C}(\theta)$ with $r(\epsilon^{\mathit{trom}}_{\mathit{on}})$ chosen via \Cref{e:energy-bound-r}
\State $\hat{X}(\theta)\gets \bigl(\Phi^{(x)}(\theta)U_r(\theta)\bigr)\,\Sigma_r(\theta)\,V_r(\theta)^\mathsf{T}$
\State $u_h^j(\theta)\gets \hat{X}(\theta)_{:,j}$ for all $j = 1,\ldots,n_t$
\State \textbf{output:} $u_h^j(\theta)$
\end{algorithmic}
\end{algorithm}

\begin{algorithm}
\caption{Offline pre-compression in space.\label{a:precompression-offline}}
\begin{algorithmic}[1]
\State \textbf{input:} tolerance $\epsilon^{\mathit{pre}}\in(0,1)$; training sets $\Theta_s=\{\theta^l\}_{l\in\mathcal{I}_s}$; $\Theta_{\mathit{pre}}=\{\theta^l\}_{l\in\mathcal{I}_{\mathit{pre}}}$
\State compute \fom\ snapshots $u_h^j(\theta^l)$ for all $\theta^l\in\Theta_{\mathit{pre}}$, $j=1,\ldots,n_t$, and assemble $S_{\mathit{pre}}$
\State $\Phi^{(x)}_{\mathit{pre}}\Sigma_{\mathit{pre}}V_{\mathit{pre}}^\mathsf{T} \gets$ compute thin SVD of $S_{\mathit{pre}}$ with $r_x(\epsilon^{\mathit{pre}})$ chosen based on \Cref{e:energy-bound-r}
\For{each $\theta^l\in\Theta_s$ and $j=1,\ldots,n_t$}
    \State $\hat{u}_h^j(\theta^l) \gets (\Phi^{(x)}_{\mathit{pre}})^\mathsf{T}u_h^j(\theta^l)$
\EndFor
\State assemble reduced $r_x\times n_t\times n_{\theta_1}\times\cdots\times n_{\theta_m}$ snapshot tensor $\hat{\mathcal{S}}$
\State \textbf{output:} pre-compression basis $\Phi^{(x)}_{\mathit{pre}}$ and compressed snapshot tensor $\hat{\mathcal{S}}$
\end{algorithmic}
\end{algorithm}

\end{appendix}

\printbibliography

@article{mang2020:integrated,
	author = {A. Mang and S. Bakas and S. Subramanian and C. Davatzikos and G. Biros},
	journal = {Annual Review of Biomedical Engineering},
	pages = {309--341},
	title = {Integrated biophysical modeling and image analysis: {A}pplication to neuro-oncology},
	volume = {22},
	year = {2020}}

@misc{tensortoolbox2025,
	author = {B. W. Bader and T. G. Kolda and others},
	howpublished = {\url{https://www.tensortoolbox.org}},
	note = {Accessed: 2026-02-23},
	title = {{Tensor Toolbox for MATLAB, Version 3.8}},
	year = {2025}}

@inproceedings{collins1999:animal,
	author = {Collins, D. L. and Zijdenbos, A. P. and Baar{\'e}, W. F. C. and Evans, A. C.},
	booktitle = {International Conference on Information Processing in Medical Imaging},
	organization = {Springer},
	pages = {210--223},
	title = {{ANIMAL+ INSECT}: {I}mproved cortical structure segmentation},
	year = {1999}}

@article{fonov2009:unbiased,
	author = {Fonov, V. S. and Evans, A. C. and McKinstry, R. C. and Almli, C. R. and Collins, D. L.},
	journal = {NeuroImage},
	pages = {S102},
	publisher = {Elsevier},
	title = {Unbiased nonlinear average age-appropriate brain templates from birth to adulthood},
	volume = {47},
	year = {2009}}

@article{fonov2011:unbiased,
	author = {Fonov, V. and Evans, A. C. and Botteron, K. and Almli, C. R. and McKinstry, R. C. and Collins, D. L. and BDCG},
	journal = {Neuroimage},
	number = {1},
	pages = {313--327},
	publisher = {Elsevier},
	title = {Unbiased average age-appropriate atlases for pediatric studies},
	volume = {54},
	year = {2011}}

@article{willcox2021:imperative,
	author = {Willcox, K. E and Ghattas, O. and Heimbach, P.},
	journal = {Nature Computational Science},
	number = {3},
	pages = {166--168},
	publisher = {Nature Publishing Group US New York},
	title = {The imperative of physics-based modeling and inverse theory in computational science},
	volume = {1},
	year = {2021}}

@article{buithanh2008:modelreduction,
	author = {Bui-Thanh, T. and Willcox, K. and Ghattas, O.},
	journal = {SIAM Journal on Scientific Computing},
	number = {6},
	pages = {3270--3288},
	title = {Model reduction for large-scale systems with polynomial nonlinearities},
	volume = {30},
	year = {2008}}

@book{hesthaven2016:certified,
	author = {Hesthaven, J. S. and Rozza, G. and Stamm, B.},
	publisher = {Springer},
	title = {Certified reduced basis methods for parametrized partial differential equations},
	volume = {590},
	year = {2016}}

@article{barrault2004:eim,
	author = {Barrault, M. and Maday, Y. and Nguyen, N. C. and Patera, A. T.},
	journal = {Comptes Rendus Mathematique},
	number = {9},
	pages = {667--672},
	title = {An empirical interpolation method: {A}pplication to efficient reduced-basis discretization of partial differential equations},
	volume = {339},
	year = {2004}}

@inproceedings{pati2020:estimating,
	author = {Pati, S. and Sharma, V. and Aslam, H. and Thakur, S. P and Akbari, H. and Mang, A. and Subramanian, S. and Biros, G. and Davatzikos, C. and Bakas, S.},
	booktitle = {International MICCAI Brainlesion Workshop},
	organization = {Springer},
	pages = {157--167},
	title = {Estimating glioblastoma biophysical growth parameters using deep learning regression},
	year = {2020}}

@article{chaturantabut2010:deim,
	author = {Chaturantabut, S. and Sorensen, D. C.},
	journal = {SIAM Journal on Scientific Computing},
	number = {5},
	pages = {2737--2764},
	title = {Nonlinear model reduction via discrete empirical interpolation},
	volume = {32},
	year = {2010}}

@article{stefanescu2014:tensorialpod,
	author = {{\c{S}}tef{\u{a}}nescu, R. and Sandu, A. and Navon, I. M.},
	journal = {International Journal for Numerical Methods in Fluids},
	number = {8},
	pages = {497--521},
	title = {Comparison of {POD} reduced order strategies for the nonlinear {2D} shallow water equations},
	volume = {76},
	year = {2014}}

@book{benner2017:mra,
	editor = {Benner, P. and Cohen, A. and Ohlberger, M. and Willcox, K.},
	publisher = {SIAM},
	title = {Model reduction and approximation: {T}heory and algorithms},
	year = {2017}}

@article{metz2024:deepgrowth,
	author = {Metz, M.-C. and Ezhov, I. and Peeken, J. C. and Buchner, J. A. and Lipkova, J. and Kofler, F. and Waldmannstetter, D. and Delbridge, C. and Diehl, C. and Bernhardt, D. and Schmidt-Graf, F. and Gempt, J. and Combs, S. E. and Zimmer, C. and Menze, B. and Wiestler, B.},
	journal = {Neuro-Oncology Advances},
	number = {1},
	pages = {vdad171},
	title = {Toward image-based personalization of glioblastoma therapy: {A} clinical and biological validation study of a novel, deep learning-driven tumor growth model},
	volume = {6},
	year = {2024}}

@article{zhang2024:pinninfiltration,
	author = {Zhang, R. Zirui and Ezhov, I. and Balcerak, M. and Zhu, A. and Wiestler, B. and Menze, B. and Lowengrub, J. S.},
	journal = {Medical Image Analysis},
	title = {Personalized Predictions of Glioblastoma Infiltration: {M}athematical Models, Physics-Informed Neural Networks and Multimodal Scans},
	year = {2024}}

@article{chen2023:tgmnet,
	author = {Chen, Q. and Ye, Q. and Zhang, W. and Li, H. and Zheng, X.},
	journal = {Engineering Applications of Artificial Intelligence},
	pages = {106867},
	title = {TGM-Nets: A deep learning framework for enhanced forecasting of tumor growth by integrating imaging and modeling},
	volume = {127},
	year = {2023}}

@article{chen2025:deep,
	author = {Chen, Q. and Li, H. and Zheng, X.},
	journal = {Engineering with Computers},
	number = {1},
	pages = {423--533},
	publisher = {Springer},
	title = {A deep neural network for operator learning enhanced by attention and gating mechanisms for long-time forecasting of tumor growth},
	volume = {41},
	year = {2025}}

@inproceedings{weidner2025:trainforwards,
	author = {Weidner, J. and Ezhov, I. and Balcerak, M. and Menze, B. and Rueckert, D. and Wiestler, B.},
	booktitle = {NeurIPS 2025 Workshop on Differentiable Simulators (DiffSys)},
	title = {Train Forwards, Optimize Backwards: Neural Surrogates for Personalized Medical Simulations},
	year = {2025}}

@article{laslo2025:diffusion,
	author = {Laslo, D. and Georgiou, E. and Linguraru, M. G. and Rauschecker, A. and Muller, S. and Jutzeler, C. R. and Bruningk, S.},
	journal = {arXiv preprint},
	title = {Mechanistic learning with guided diffusion models to predict spatio-temporal brain tumor growth},
	year = {2025}}

@article{liu2025:m4rl,
	author = {Liu, Z. and Zhang, J. and Hong, L. and Nie, Q. and Sun, X.},
	journal = {Science Advances},
	number = {32},
	pages = {eadv3316},
	title = {Multiscale mathematical model-informed reinforcement learning optimizes combination treatment scheduling in glioblastoma evolution},
	volume = {11},
	year = {2025}}

@article{benner2015:survey,
	author = {Benner, P. and Gugercin, S. and Willcox, K.},
	journal = {SIAM Review},
	number = {4},
	pages = {483--531},
	publisher = {SIAM},
	title = {A survey of projection-based model reduction methods for parametric dynamical systems},
	volume = {57},
	year = {2015}}

@article{dolgov2012:fast,
	author = {Dolgov, S. V. and Khoromskij, B. N. and Oseledets, I. V.},
	journal = {SIAM Journal on Scientific Computing},
	number = {6},
	pages = {A3016--A3038},
	publisher = {SIAM},
	title = {Fast solution of parabolic problems in the tensor train/quantized tensor train format with initial application to the Fokker--Planck equation},
	volume = {34},
	year = {2012}}

@article{manzini2025:low,
	author = {Manzini, G. and Sorgente, T.},
	journal = {Mathematics and Computers in Simulation},
	publisher = {Elsevier},
	title = {The low-rank tensor-train finite difference method for three-dimensional parabolic equations},
	year = {2025}}

@article{chaturantabut2010:nonlinear,
	author = {Chaturantabut, S. and Sorensen, D. C.},
	journal = {SIAM Journal on Scientific Computing},
	number = {5},
	pages = {2737--2764},
	title = {Nonlinear model reduction via discrete empirical interpolation},
	volume = {32},
	year = {2010}}

@article{halko2011:finding,
	author = {Halko, N. and Martinsson, P.-G. and Tropp, J. A.},
	journal = {SIAM Review},
	number = {2},
	pages = {217--288},
	title = {Finding structure with randomness: {P}robabilistic algorithms for constructing approximate matrix decompositions},
	volume = {53},
	year = {2011}}

@book{mang2014:methoden,
	author = {Mang, A.},
	publisher = {Springer-Verlag},
	title = {Methoden zur numerischen {S}imulation der {P}rogression von {G}liomen: {M}odellentwicklung, {N}umerik und {P}arameteridentifikation},
	year = {2014}}

@article{oden2010:general,
	author = {J. T. Oden and A. Hawkins and S. Prudhomme},
	journal = {Mathematical Models and Methods in Applied Sciences},
	number = {3},
	pages = {477--517},
	title = {General diffuse-interface theories and an approach to predictive tumor growth modeling},
	volume = {20},
	year = {2010}}

@article{bakas2018:identifying,
	archiveprefix = {arXiv},
	author = {S. Bakas and M. Reyes and A. Jakab and S. Bauer and M. Rempfler and others},
	eprint = {1811.02629},
	journal = {arXiv e-prints},
	primaryclass = {cs.CV},
	title = {Identifying the best machine learning algorithms for brain tumor segmentation, progression assessment, and overall survival prediction in the {BRATS} challenge},
	year = {2018}}

@article{bakas2015:glistrboost,
	author = {S. Bakas and Z. Zeng and A. Sotiras and S. Rathore and H. Akbari and B. Gaonkar and M. Rozycki and S. Pati and C. Davatzikos},
	journal = {Brain Lesion},
	pages = {144--155},
	title = {{GLISTR}boost: {C}ombining multimodal {MRI} segmentation, registration, and biophysical tumor growth modeling with gradient boosting machines for glioma segmentation},
	volume = {9556},
	year = {2015}}

@article{oden2016:toward,
	author = {Oden, J. T. and Lima, E. A. B. F. and Almeida, R. C. and Feng, Y. and Rylander, M. N. and Fuentes, D. and Faghihi, D. and Rahman, M. M. and DeWitt, M. and Gadde, M. and Zhou, J. C.},
	journal = {Archives of Computational Methods in Engineering},
	number = {4},
	pages = {735--779},
	publisher = {Springer},
	title = {Toward predictive multiscale modeling of vascular tumor growth},
	volume = {23},
	year = {2016}}

@article{yankeelov2013:clinically,
	author = {Yankeelov, T. E. and Atuegwu, N. and Hormuth, D. and Weis, J. A. and Barnes, S. L. and Miga, M. I. and Rericha, E. C. and Quaranta, V.},
	journal = {Science Translational Medicine},
	number = {187},
	pages = {187ps9},
	title = {Clinically relevant modeling of tumor growth and treatment response},
	volume = {5},
	year = {2013}}

@article{jbabdi2005:simulation,
	author = {Jbabdi, S. and Mandonnet, E. and Duffau, H. and Capelle, L. and Swanson, K. R. and P{\'e}l{\'e}grini-Issac, M. and Guillevin, R. and Benali, H.},
	journal = {Magnetic Resonance in Medicine},
	number = {3},
	pages = {616--624},
	title = {Simulation of anisotropic growth of low-grade gliomas using diffusion tensor imaging},
	volume = {54},
	year = {2005}}

@article{oden2013:selection,
	author = {J. T. Oden and E. E. Prudencio and A. Hawkins-Daarud},
	journal = {Mathematical Models and Methods in Applied Sciences},
	number = {7},
	pages = {1309--1338},
	title = {Selection and assessment of phenomenological models of tumor growth},
	volume = {23},
	year = {2013}}

@article{lima2017:selection,
	author = {E. Lima and J. T. Oden and B. Wohlmuth and A. Shahmoradi and D. A. Hormuth and T. E. Yankeelov and L. Scarabosio and T. Horger},
	journal = {Computer Methods in Applied Mechanics and Engineering},
	pages = {227--308},
	title = {Selection and validation of predictive models of radiation effects on tumor growth based on noninvasive imaging data},
	volume = {327},
	year = {2017}}

@article{jarrett2018:mathematical,
	author = {A. M. Jarrett and E. Lima and D. A. Hormuth and M. T. McKenna and X. Feng and D. A. Ekrut and A. C. M. Resende and A. Brok and T. E. Yankeelov},
	journal = {Expert Review of Anticancer Therapy},
	pages = {1271--1286},
	title = {Mathematical models of tumor cell proliferation: {A} review of the literature},
	volume = {18},
	year = {2018}}

@article{hormuth2017:mechanically,
	author = {D. A. Hormuth and J. A. Weis and S. L. Barnes and M. I. Miga and E. C. Rechricha and V. Quaranta and T. E. Yankeelov},
	journal = {Journal of the Royal Socieity Interface},
	number = {20161010},
	title = {A mechanically coupled reaction-diffusion model that incorporates intra-tumoural heterogeneity to predict in vivo glioma growth},
	volume = {14},
	year = {2017}}

@article{hawkinsdaarud2012:numerical,
	author = {A. Hawkins-Daarud and K. G. Zee and J. T. Oden},
	journal = {International Journal for Numerical Methods in Biomedical Engineering},
	number = {1},
	pages = {3--24},
	title = {Numerical simulation of a thermodynamically consistent four-species tumor growth model},
	volume = {28},
	year = {2012}}

@article{swanson2011:quantifying,
	author = {K. R. Swanson and R. C. Rockne and J. Claridge and M. A. Chaplain and E. C. Alvord and A. R. Anderson},
	journal = {Cancer Research},
	number = {24},
	pages = {7366--7375},
	title = {Quantifying the role of angiogenesis in malignant progression of gliomas: {I}n silico modeling integrates imaging and histology},
	volume = {71},
	year = {2011}}

@article{swanson2003:virtual,
	author = {K. R. Swanson and C. Bridge and J. D. Murray and E. C. Alvord},
	journal = {Journal of the Neurological Sciences},
	pages = {1-10},
	title = {Virtual and real brain tumors: using mathematical modeling to quantify glioma growth and invasion},
	volume = 216,
	year = 2003}

@article{kilmer2021:tensor,
	author = {Kilmer, M. E. and Horesh, L. and Avron, H. and Newman, E.},
	journal = {Proceedings of the National Academy of Sciences},
	number = {28},
	pages = {e2015851118},
	publisher = {National Academy of Sciences},
	title = {Tensor-tensor algebra for optimal representation and compression of multiway data},
	volume = {118},
	year = {2021}}

@article{tucker1966:some,
	author = {Tucker, L. R.},
	journal = {Psychometrika},
	number = {3},
	pages = {279--311},
	title = {Some mathematical notes on three-mode factor analysis},
	volume = {31},
	year = {1966}}

@article{de2000:multilinear,
	author = {De Lathauwer, L. and De Moor, B. and Vandewalle, J.},
	journal = {SIAM journal on Matrix Analysis and Applications},
	number = {4},
	pages = {1253--1278},
	title = {A multilinear singular value decomposition},
	volume = {21},
	year = {2000}}

@article{ghafouri2025:inverse,
	author = {Ghafouri, A. and Biros, G.},
	journal = {International Journal for Numerical Methods in Biomedical Engineering},
	number = {7},
	pages = {e70057},
	publisher = {Wiley Online Library},
	title = {Inverse problem regularization for {3D} multi-species tumor growth models},
	volume = {41},
	year = {2025}}

@article{gooya2012:glistr,
	author = {Gooya, A. and Pohl, K. M. and Bilello, M. and Cirillo, L. and Biros, G. and Melhem, E. R. and Davatzikos, C.},
	journal = {IEEE Transactions on Medical Imaging},
	number = {10},
	pages = {1941--1954},
	publisher = {IEEE},
	title = {{GLISTR}: {G}lioma image segmentation and registration},
	volume = {31},
	year = {2012}}

@article{yadav2025:understanding,
	author = {Yadav, V. S. and Ranwan, N. and Chamakuri, N.},
	journal = {Computers \& Mathematics with Applications},
	pages = {55--70},
	publisher = {Elsevier},
	title = {Understanding avascular tumor growth and drug interactions through numerical analysis: {A} finite element method approach},
	volume = {181},
	year = {2025}}

@article{krawczyk2020:optimal,
	author = {Krawczyk, A. and Nowakowski, A.},
	journal = {Computers \& Mathematics With Applications},
	number = {5},
	pages = {778--789},
	publisher = {Elsevier},
	title = {Optimal control of inhibiting tumor, new approach to sufficient $\varepsilon$-optimality and numerical computation},
	volume = {80},
	year = {2020}}

@article{sabir2017:mathematical,
	author = {Sabir, M. and Shah, A. and Muhammad, W. and Ali, I. and Bastian, P.},
	journal = {Computers \& Mathematics with Applications},
	number = {12},
	pages = {3250--3259},
	publisher = {Elsevier},
	title = {A mathematical model of tumor hypoxia targeting in cancer treatment and its numerical simulation},
	volume = {74},
	year = {2017}}

@article{ozuugurlu2015:note,
	author = {{\"O}zu{\u{g}}urlu, E},
	journal = {Computers \& Mathematics with Applications},
	number = {12},
	pages = {1504--1517},
	publisher = {Elsevier},
	title = {A note on the numerical approach for the reaction--diffusion problem to model the density of the tumor growth dynamics},
	volume = {69},
	year = {2015}}

@article{knopoff2013:adjoint,
	author = {Knopoff, D. A. and Fern{\'a}ndez, D. R. and Torres, G. A. and Turner, C. V.},
	journal = {Computers \& Mathematics with Applications},
	number = {6},
	pages = {1104--1119},
	publisher = {Elsevier},
	title = {Adjoint method for a tumor growth {PDE}-constrained optimization problem},
	volume = {66},
	year = {2013}}

@article{hogea2008:image,
	author = {Hogea, C. and Davatzikos, C. and Biros, G.},
	journal = {Journal of Mathematical Biology},
	number = {6},
	pages = {793--825},
	publisher = {Springer},
	title = {An image-driven parameter estimation problem for a reaction--diffusion glioma growth model with mass effects},
	volume = {56},
	year = {2008}}

@article{scheufele2020:imagedriven,
	author = {K. Scheufele and S. Subramanian and A. Mang and G. Biros and M. Mehl},
	journal = {SIAM Journal on Scientific Computing},
	number = {3},
	pages = {B549--B580},
	title = {Image-driven biophysical tumor growth model calibration},
	volume = {42},
	year = {2020}}

@article{scheufele2019:coupling,
	author = {K. Scheufele and A. Mang and A. Gholami and C. Davatzikos and G. Biros and M. Mehl},
	journal = {Computer Methods in Applied Mechanics and Engineering},
	pages = {533--567},
	title = {Coupling brain-tumor biophysical models and diffeomorphic image registration},
	volume = {347},
	year = {2019}}

@inproceedings{gholami2017:framework,
	author = {Gholami, A. and Mang, A. and Scheufele, K. and Davatzikos, C. and Mehl, M. and Biros, G.},
	booktitle = {Proc ACM/IEEE Conference on Supercomputing},
	number = {19},
	pages = {19:1--19:13},
	title = {A framework for scalable biophysics-based image analysis},
	year = {2017}}

@article{mang2018:pdeconst,
	author = {Mang, A. and Gholami, A. and Davatzikos, C. and Biros, G.},
	journal = {Optimization and Engineering},
	number = {3},
	pages = {765--812},
	title = {{PDE}-constrained optimization in medical image analysis},
	volume = {19},
	year = {2018}}

@article{mang2012:biophys,
	author = {Mang, A. and Toma, A. and Sch\"utz, T. A. and Becker, S. and Eckey, T. and Mohr, C. and Petersen, D. and Buzug, T. M.},
	journal = {Medical Physics},
	number = {7},
	pages = {4444--4459},
	title = {Biophysical modeling of brain tumor progression: from unconditionally stable explicit time integration to an inverse problem with parabolic {PDE} constraints for model calibration},
	volume = {39},
	year = {2012}}

@article{gholami2016:inverse,
	author = {Gholami, A. and Mang, A. and Biros, G.},
	journal = {Journal of Mathematical Biology},
	number = {1},
	pages = {409--433},
	publisher = {Springer},
	title = {An inverse problem formulation for parameter estimation of a reaction--diffusion model of low grade gliomas},
	volume = {72},
	year = {2016}}

@article{biros2023:inverse,
	author = {Biros, G. and Mang, A. and Menze, B. H. and Schulte, M.},
	journal = {Dagstuhl Reports},
	number = {1},
	pages = {36--67},
	publisher = {Schloss Dagstuhl--Leibniz-Zentrum f{\"u}r Informatik},
	title = {Inverse biophysical modeling and machine learning in personalized oncology ({D}agstuhl {S}eminar 23022)},
	volume = {13},
	year = {2023}}

@article{subramanian2019:simulation,
	author = {Subramanian, S. and Gholami, A. and Biros, G.},
	journal = {Journal of Mathematical Biology},
	pages = {941--967},
	publisher = {Springer},
	title = {Simulation of glioblastoma growth using a {3D} multispecies tumor model with mass effect},
	volume = {79},
	year = {2019}}

@article{saut2014:multilayer,
	author = {Saut, O. and Lagaert, J.-B. and Colin, T. and Fathallah-Shaykh, H. M.},
	journal = {Bulletin of Mathematical Biology},
	pages = {2306--2333},
	publisher = {Springer},
	title = {A multilayer grow-or-go model for {GBM}: {E}ffects of invasive cells and anti-angiogenesis on growth},
	volume = {76},
	year = {2014}}

@book{hundsdorfer2003:numerical,
	author = {Hundsdorfer, W. and Verwer, J. G.},
	publisher = {Springer},
	title = {Numerical solution of time-dependent advection--diffusion--reaction equations},
	year = {2003}}

@article{strang1968:difference,
	author = {Strang, G.},
	journal = {SIAM Journal on Numerical Analysis},
	number = {3},
	pages = {506--517},
	publisher = {SIAM},
	title = {On the construction and comparison of difference schemes},
	volume = {5},
	year = {1968}}

@article{eymard2000finite,
	author = {Eymard, R. and Gallou{\"e}t, T. and Herbin, R.},
	journal = {Handbook of Numerical Analysis},
	pages = {713--1020},
	publisher = {Elsevier},
	title = {The finite volume method},
	volume = {7},
	year = {2000}}

@book{leveque2002:fvm,
	author = {LeVeque, R. J.},
	publisher = {Cambridge University Press},
	title = {Finite Volume Methods for Hyperbolic Problems},
	year = {2002}}

@techreport{kolda2017:tensor,
	author = {Kolda, T. and Bader, B. W. and Acar A., Evrim N. M. N. and Dunlavy, D. and Bassett, R. and Battaglino, C. J. and Plantenga, T. and Chi, E. and Hansen, S.},
	institution = {Sandia National Laboratories (SNL-NM), Albuquerque, NM (United States)},
	title = {Tensor Toolbox for {MATLAB}},
	year = {2017}}

@article{bader2008:efficient,
	author = {Bader, B. W. and Kolda, T. G.},
	journal = {SIAM Journal on Scientific Computing},
	number = {1},
	pages = {205--231},
	publisher = {SIAM},
	title = {Efficient {MATLAB} computations with sparse and factored tensors},
	volume = {30},
	year = {2008}}

@article{li2020:fourier,
	author = {Li, Z. and Kovachki, N. and Azizzadenesheli, K. and Liu, B. and Bhattacharya, K. and Stuart, A. and Anandkumar, A.},
	journal = {arXiv preprint arXiv:2010.08895},
	title = {Fourier neural operator for parametric partial differential equations},
	year = {2020}}

@article{lu2021:learning,
	author = {Lu, L. and Jin, P. and Pang, G. and Zhang, Z. and Karniadakis, G. E.},
	journal = {Nature Machine Intelligence},
	number = {3},
	pages = {218--229},
	publisher = {Nature Publishing Group UK London},
	title = {Learning nonlinear operators via {DeepONet} based on the universal approximation theorem of operators},
	volume = {3},
	year = {2021}}

@article{raissi2019:physics,
	author = {Raissi, M. and Perdikaris, P. and Karniadakis, G. E.},
	journal = {Journal of Computational physics},
	pages = {686--707},
	publisher = {Elsevier},
	title = {Physics-informed neural networks: {A} deep learning framework for solving forward and inverse problems involving nonlinear partial differential equations},
	volume = {378},
	year = {2019}}

@book{antoulas2005:approximation,
	author = {Antoulas, A. C.},
	publisher = {SIAM},
	title = {Approximation of large-scale dynamical systems},
	year = {2005}}

@article{rozza2008:reduced,
	author = {Rozza, G. and Huynh, D. B. P. and Patera, A. T.},
	journal = {Archives of Computational Methods in Engineering},
	number = {3},
	pages = {229--275},
	publisher = {Springer},
	title = {Reduced basis approximation and a posteriori error estimation for affinely parametrized elliptic coercive partial differential equations: {A}pplication to transport and continuum mechanics},
	volume = {15},
	year = {2008}}

@article{peherstorfer2018:survey,
	author = {Peherstorfer, B. and Willcox, K. and Gunzburger, M.},
	journal = {SIAM Review},
	number = {3},
	pages = {550--591},
	publisher = {SIAM},
	title = {Survey of multifidelity methods in uncertainty propagation, inference, and optimization},
	volume = {60},
	year = {2018}}

@article{ghattas2021:learning,
	author = {Ghattas, O. and Willcox, K.},
	journal = {Acta Numerica},
	pages = {445--554},
	publisher = {Cambridge University Press},
	title = {Learning physics-based models from data: {P}erspectives from inverse problems and model reduction},
	volume = {30},
	year = {2021}}

@article{kunisch2002:galerkin,
	author = {Kunisch, K. and Volkwein, S.},
	journal = {SIAM Journal on Numerical Analysis},
	number = {2},
	pages = {492--515},
	publisher = {SIAM},
	title = {Galerkin proper orthogonal decomposition methods for a general equation in fluid dynamics},
	volume = {40},
	year = {2002}}

@book{quarteroni2015:reduced,
	author = {Quarteroni, A. and Manzoni, A. and Negri, F.},
	publisher = {Springer},
	title = {Reduced basis methods for partial differential equations: {A}n introduction},
	year = {2015}}

@article{mamonov2025:priori,
	author = {Mamonov, A. V. and Olshanskii, M. A.},
	journal = {SIAM Journal on Numerical Analysis},
	number = {1},
	pages = {239--261},
	publisher = {SIAM},
	title = {A priori analysis of a tensor {ROM} for parameter dependent parabolic problems},
	volume = {63},
	year = {2025}}

@article{mamonov2024slice,
	author = {Mamonov, Alexander V and Olshanskii, Maxim A},
	journal = {arXiv e-prints},
	pages = {arXiv--2411},
	title = {Slice sampling tensor completion for model order reduction of parametric dynamical systems},
	year = {2024}}

@article{budzinskiy2025low,
	author = {Budzinskiy, S. and Kazeev, V. and Olshanskii, M.},
	journal = {arXiv preprint arXiv:2511.22650},
	title = {Low-rank cross approximation of function-valued tensors for reduced-order modeling of parametric {PDE}s},
	year = {2025}}

@article{mamonov2022:interpolatory,
	author = {Mamonov, A. V. and Olshanskii, M. A.},
	journal = {Computer Methods in Applied Mechanics and Engineering},
	pages = {115122},
	publisher = {Elsevier},
	title = {Interpolatory tensorial reduced order models for parametric dynamical systems},
	volume = {397},
	year = {2022}}

@article{mueller2025:tensor,
	author = {Mueller, N. and Zhao, Y. and Badia, S. and Cui, T.},
	journal = {Journal of Computational and Applied Mathematics},
	pages = {116790},
	publisher = {Elsevier},
	title = {A tensor-train reduced basis solver for parameterized partial differential equations on Cartesian grids},
	year = {2025}}

@article{olshanskii2025:approximating,
	author = {Olshanskii, M. A. and Rebholz, L. G.},
	journal = {Journal of Computational Physics},
	pages = {113728},
	publisher = {Elsevier},
	title = {Approximating a branch of solutions to the {N}avier--{S}tokes equations by reduced-order modeling},
	volume = {524},
	year = {2025}}

@article{mamonov2024:tensorial,
	author = {Mamonov, A. V. and Olshanskii, M. A.},
	journal = {SIAM Journal on Scientific Computing},
	number = {3},
	pages = {A1850--A1878},
	publisher = {SIAM},
	title = {Tensorial parametric model order reduction of nonlinear dynamical systems},
	volume = {46},
	year = {2024}}

@article{mizan2025:parametric,
	author = {Md {Rezwan Bin Mizan} and Maxim Olshanskii and Ilya Timofeyev},
	journal = {Computers \& Fluids},
	pages = {107005},
	title = {A parametric tensor ROM for the shallow water dam break problem},
	volume = {309},
	year = {2026}}

@article{scheufele2020:fully,
	author = {Scheufele, K. and Subramanian, S. and Biros, G.},
	journal = {IEEE Transactions on Medical Imaging},
	number = {1},
	pages = {193--204},
	publisher = {IEEE},
	title = {Fully automatic calibration of tumor-growth models using a single {mpMRI} scan},
	volume = {40},
	year = {2020}}

@article{oseledets2011:tensor_train,
	author = {Oseledets, I. V.},
	journal = {SIAM Journal on Scientific Computing},
	number = {5},
	pages = {2295--2317},
	publisher = {SIAM},
	title = {Tensor-{T}rain Decomposition},
	volume = {33},
	year = {2011}}

@article{zhang2025:personalized,
	author = {Zhang, R. Z. and Ezhov, I. and Balcerak, M. and Zhu, A. and Wiestler, B. and Menze, B. and Lowengrub, J. S.},
	journal = {Medical Image Analysis},
	pages = {103423},
	publisher = {Elsevier},
	title = {Personalized predictions of {G}lioblastoma infiltration: {M}athematical models, physics-informed neural networks and multimodal scans},
	volume = {101},
	year = {2025}}

\end{document}